\documentclass[a4paper,12pt]{article}     

\usepackage{graphicx}		  
\usepackage{amsmath,amssymb}	

\newtheorem{lemma}{Lemma}[section]
\newtheorem{theorem}[lemma]{Theorem}
\newtheorem{proposition}[lemma]{Proposition}

\newtheorem{assumption}{Assumption}

\newcommand\LABEL[1]{\label{#1}}

\def\authorfont{\footnotesize}

\def\ccode#1{\par
\vspace*{8pt}
{\authorfont{\leftskip18pt\rightskip\leftskip
\noindent #1\par}}\par}

\newenvironment{Proof}{
\hspace*{-9mm}
{ \it Proof.}}
{\hfill {$\square$}\vspace{1.5em}}

\begin{document}

\begin{center}{
{\Large 
 Properties of minimal charts and
 their applications VIII: 
 charts of type $(7)$}
\vspace{10pt}
\\ 
Teruo NAGASE and Akiko SHIMA\footnote{The second author is supported by JSPS KAKENHI Grant Number 21K03255.}
}
\end{center}


\begin{abstract}
Let $\Gamma$ be a chart,
and we denote by $\Gamma_m$
the union of all the edges of label $m$.
A chart $\Gamma$ is of type $(7)$
if there exists a label $m$
such that 
$w(\Gamma)=7$,
$w(\Gamma_m\cap\Gamma_{m+1})=7$
where 
$w(G)$ is the number of white vertices in $G$.
In this paper, we prove that there is 
no minimal chart of 
type $(7)$.
\end{abstract}

%
%
%
%

\ccode{2020 Mathematics Subject Classification. Primary 05C10; Secondary 57M15.}
\ccode{ {\it Key Words and Phrases}. surface link, chart, white vertex. }

\setcounter{section}{0}
\section{Introduction}


Charts are oriented labeled graphs in a disk (see  \cite{KnottedSurfaces},\cite{BraidBook}, and see Section~\ref{s:Prel}  for the precise definition of charts).
From a chart, we can construct an oriented closed surface 
embedded in 4-space ${\Bbb R}^4$ 
 (see \cite[Chapter 14, Chapter 18 and Chapter 23]{BraidBook}). 
A C-move 
is a local modification between two charts
in a disk (see Section~\ref{s:Prel} for C-moves).
A C-move between two charts induces 
an ambient isotopy between oriented closed surfaces 
corresponding to the two charts.

We will work in the PL category or smooth category. All submanifolds are assumed to be locally flat.
In \cite{ONS},
we showed that there is no minimal chart with exactly five vertices
 (see Section~\ref{s:Prel} for the precise definition of minimal charts). 
Hasegawa proved that there exists a minimal chart with exactly
six white vertices \cite{H1}. 
This chart represents a 2-twist spun trefoil.
In \cite{INS} and \cite{NST},
we investigated minimal charts with exactly four white vertices.
In this paper, 
we investigate properties of minimal charts and
need to prove that
there is no minimal chart with exactly seven white vertices
(see \cite{ChartApp1},\cite{ChartAppII},
\cite{ChartAppIII},\cite{ChartAppIV},
\cite{ChartAppV}, \cite{ChartAppVI},
\cite{ChartAppVII},\cite{ChartAppIX}).

Let $\Gamma$ be a chart.
For each label $m$, we denote by $\Gamma_m$
the union of all the edges of label $m$.

Now we define a type of a chart:
Let $\Gamma$ be a chart with at least one white vertex, 
and $n_1,n_2,\dots,n_k$ integers.
The chart $\Gamma$ is of {\it type $(n_1,n_2,\dots,n_k)$} if there exists a label $m$ of $\Gamma$ satisfying the following three conditions:
\begin{enumerate}
\item[(i)] For each $i=1,2,\dots, k$, 
the chart $\Gamma$ contains exactly $n_{i}$ white vertices in $\Gamma_{m+i-1}\cap \Gamma_{m+i}$.
\item[(ii)] If $i<0$ or $i>k$, then $\Gamma_{m+i}$ does not contain any white vertices.
\item[(iii)] Both of the two subgraphs $\Gamma_m$ and $\Gamma_{m+k}$ contain at least one white vertex.
\end{enumerate}
If we want to emphasize the label $m$,
then we say that $\Gamma$ is of {\it type $(m;n_1,n_2,\dots,n_k)$}. 
Note that $n_1\ge1$ and $n_k\ge1$ by the condition (iii).

We proved in \cite[Theorem 1.1]{ChartAppII} that
if there exists a minimal $n$-chart $\Gamma$ with exactly seven white vertices,
then $\Gamma$ is a chart of 
type $(7),(5,2),(4,3),(3,2,2)$ or $(2,3,2)$ 
(if necessary we change the label
$i$ by $n-i$ for all label $i$).
In \cite{ChartAppV},
we showed that
there is no minimal chart of type $(3,2,2)$.
In \cite{ChartAppVI} and \cite{ChartAppVII},
there is no minimal chart of type $(2,3,2)$.
In this paper we shall show the following:

\begin{theorem}
\LABEL{MainTheorem} 
There is 
no minimal chart of 
type $(7)$.
\end{theorem}

In the future paper \cite{ChartAppIX},
we shall show there is no minimal chart of type
$(5,2),(4,3)$.
Therefore we shall show that
there is no minimal chart with exactly seven white vertices.

The paper is organized as follows.
In Section~\ref{s:Prel},
we define charts and minimal charts.
In Section~\ref{s:2-angledDisks},
we investigate a disk called a 2-angled disk
whose boundary contains exactly two white vertices
and consists of edges of label $m$.
In Section~\ref{s:3-angledDisks}, 
we review a lemma about a disk called a 3-angled disk
whose boundary contains exactly three white vertices
and consists of edges of label $m$.
In Section~\ref{s:TypeGammaM}, 
we investigate the graph $\Gamma_m$
for a minimal chart $\Gamma$ of type $(m;7)$.
In Section~\ref{s:Rings},
we shall show that 
if $\Gamma$ is a minimal chart of type $(m;7)$,
then there exist two connected components $G_1$ and $G_2$
of $\Gamma_m$ (resp. $\Gamma_{m+1}$) such that
$G_1$ contains five white vertices and
$G_2$ contains two white vertices.
In Section~\ref{s:FiveWhiteVertices},
we investigate
a connected component of $\Gamma_m$
with exactly five white vertices 
for a minimal chart $\Gamma$.
In Section~\ref{s:TypeATypeD},
we shall show that
neither $\Gamma_{m}$ nor $\Gamma_{m+1}$ contains
the graphs as shown in Fig.~\ref{Fig12}(a) and (d)
for any minimal chart $\Gamma$ of type $(m;7)$.
In Section~\ref{s:TypeF},
we review IO-Calculation(a property of numbers of 
inward arcs of label $k$ 
and outward arcs of label $k$ in a closed domain $F$
with $\partial F\subset\Gamma_{k-1}\cup\Gamma_k\cup\Gamma_{k+1}$
for some label $k$).
Moreover we shall show that
neither $\Gamma_{m}$ nor $\Gamma_{m+1}$ contains
the graph as shown in Fig.~\ref{Fig12}(f)
for any minimal chart $\Gamma$ of type $(m;7)$.
From Section~\ref{s:TypeI} to Section~\ref{s:TypeC},
we shall show that
neither $\Gamma_{m}$ nor $\Gamma_{m+1}$ contains
the graphs as shown in Fig.~\ref{Fig12}(i),(h),(e),(b),(c)
for any minimal chart $\Gamma$ of type $(m;7)$.
In Section~\ref{s:Main},
we shall prove Theorem~\ref{MainTheorem}.


\section{Preliminaries}
\LABEL{s:Prel}

In this section, 
we introduce 
the definition of charts and its related words.

Let $n$ be a positive integer.
An $n$-{\it chart}  
(a braid chart of degree $n$ \cite{KnottedSurfaces}
or a surface braid chart of degree $n$ \cite{BraidBook}) 
is 
an oriented labeled graph in the interior of a disk,
which may be empty 
or
have closed edges without vertices
satisfying the following four conditions
(see Fig.~\ref{Fig01}):
\begin{enumerate}
\item[(i)] 
Every vertex has degree $1$, $4$, or $6$.
\item[(ii)] 
The labels of edges are 
in $\{1,2,\dots,n-1\}$.
\item[(iii)]
In a small neighborhood of
each vertex of degree $6$,
there are six short arcs,
three consecutive arcs are
oriented inward 
and
the other three are outward,
and
these six are labeled $i$ and $i+1$
alternately for some $i$,
where the orientation and label of
each arc are inherited from
the edge containing the arc.
\item[(iv)]
For each vertex of degree $4$,
diagonal edges have the same label
and
are oriented coherently,
and the labels $i$ and $j$ of
the diagonals satisfy $|i-j|>1$.
\end{enumerate}
We call a vertex of degree $1$ a {\it black vertex},
a vertex of degree $4$ a {\it crossing}, and 
a vertex of degree $6$ a {\it white vertex}
respectively.

Among six short arcs
in a small neighborhood of
a white vertex,
a central arc of each three consecutive arcs
oriented inward (resp. outward) 
is called a   
{\it middle arc} at the white vertex
(see Fig.~\ref{Fig01}(c)).
For each white vertex $v$, 
there are two middle arcs at $v$ 
in a small neighborhood of $v$.
An edge is said to be {\it middle at} a white vertex $v$ if it contains a middle arc at $v$.

Let $e$ be an edge connecting $v_1$ and $v_2$.
If $e$ is oriented from $v_1$ to $v_2$,
then we say that 
$e$ is oriented {\it outward at $v_1$}
and {\it inward at $v_2$}


\begin{figure}[htb]
\begin{center}
\includegraphics{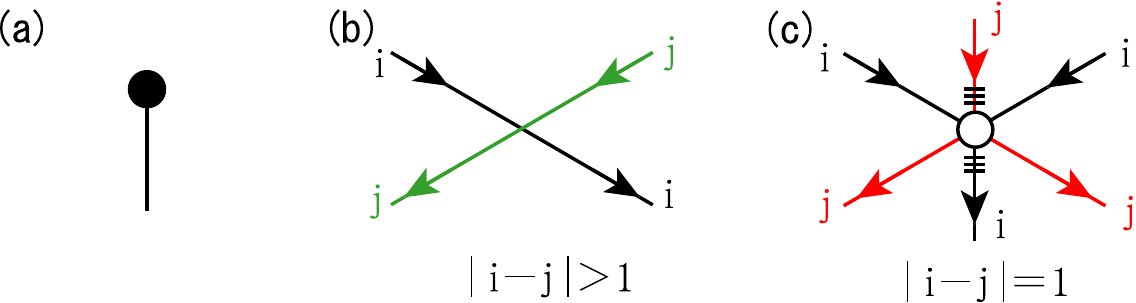}
\end{center}
\caption{ \LABEL{Fig01} (a) A black vertex. (b) A crossing. (c) A white vertex. 
Each arc with three transversal short arcs is a middle arc at the white vertex. }
\end{figure}

Now {\it C-moves} are local modifications 
of charts as shown in Fig.~\ref{Fig02}
(cf. \cite{KnottedSurfaces}, 
\cite{BraidBook} and \cite{Tanaka}).
Two charts are said to be {\it C-move equivalent}  if there exists
a finite sequence of C-moves 
which modifies one of the two charts 
to the other.

\begin{figure}
\begin{center}
\includegraphics{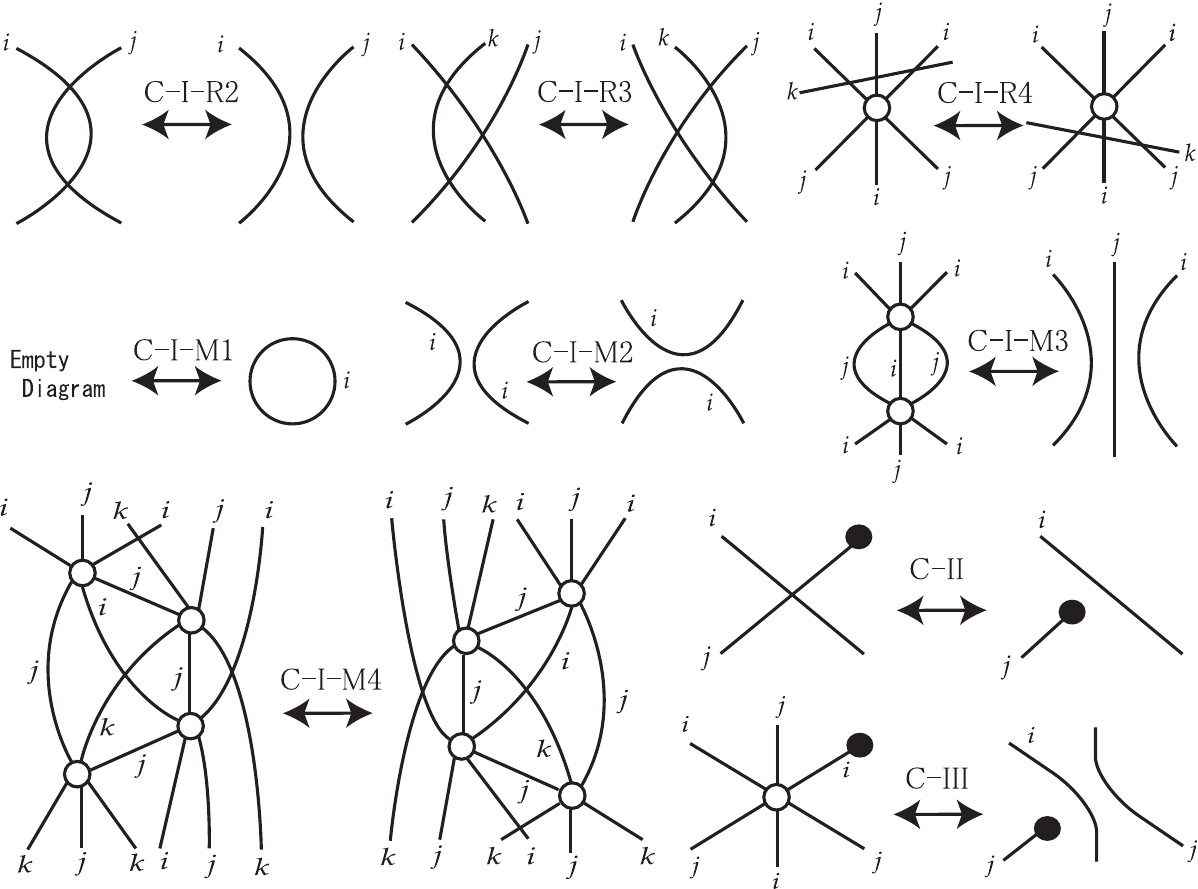}
\end{center}
\caption{ \LABEL{Fig02} For the C-III move, 
the edge with the black vertex is not middle at
a white vertex in the left figure. }
\end{figure}

An edge in a chart is called 
a {\it free edge}
if it has
two black vertices.

For each chart $\Gamma$,
let $w(\Gamma)$ and $f(\Gamma)$ be the number of white vertices, and the number of free edges respectively.
The pair $(w(\Gamma), -f(\Gamma))$ is called a {\it complexity} of the chart (see \cite{BraidThree}).
A chart $\Gamma$ is called a {\it minimal chart} if its complexity is minimal among the charts C-move equivalent to the chart $\Gamma$ with respect to the lexicographic order of pairs of integers.

We showed the difference of a chart in a disk and in a 2-sphere (see \cite[Lemma 2.1]{ChartApp1}).
This lemma follows from that there exists a natural one-to-one correspondence between $\{$charts in $S^2\}/$C-moves and $\{$charts in $D^2\}/$C-moves, conjugations
(\cite[Chapter 23 and Chapter 25]{BraidBook}).
To make the argument simple, we assume that 
the charts lie on the 2-sphere instead of the disk.
\begin{assumption}
In this paper,
all charts are contained in the $2$-sphere $S^2$.
\end{assumption}
We have the special point in the 2-sphere $S^2$, called the point at infinity,
 denoted by $\infty$.
In this paper, all charts are contained in a disk such that the disk 
does not contain the point at infinity $\infty$.

An edge in a chart is called 
a {\it terminal edge}
if it has
a white vertex and a black vertex.

Let $\Gamma$ be a chart,
and $m$ a label of $\Gamma$. 
A {\it hoop} is a closed edge of $\Gamma$ without vertices 
(hence without crossings, neither).
A {\it ring} is a simple closed curve in $\Gamma_m$ containing a crossing but not containing any white vertices.
A hoop is said to be {\it simple} 
if one of the two complementary domains
of the hoop
does not contain any white vertices.

We can assume that
all minimal charts $\Gamma$
satisfy the following four conditions 
(see \cite{ChartApp1},\cite{ChartAppII},\cite{ChartAppIII},
\cite{StI}):

\begin{assumption}
\LABEL{AssumeTerminal}
If an edge of $\Gamma$
contains a black vertex,
then the edge is a free edge 
or a terminal edge.
Moreover 
any terminal edge contains a middle arc.
\end{assumption}

\begin{assumption}
\LABEL{NoSimpleHoop}
All free edges and simple hoops in $\Gamma$ 
are moved into a small neighborhood $U_\infty$ 
of the point at infinity $\infty$. 
Hence
we assume that 
$\Gamma$ does not contain free edges
nor simple hoops, 
otherwise mentioned. 
\end{assumption}

\begin{assumption}
\LABEL{Ring}
Each complementary domain of
any ring and hoop must contain 
at least one white vertex. 
\end{assumption}

\begin{assumption}
\LABEL{Infinity}
The point at infinity $\infty$ is moved in any complementary domain of $\Gamma$.
\end{assumption}

In this paper
for a set $X$ in a space
we denote 
the interior of $X$,
the boundary of $X$ and
the closure of $X$
by Int$X$, $\partial X$
and $Cl(X)$
respectively.

Let $\alpha$ be a simple arc or an edge,
and $p,q$ the endpoints of $\alpha$.
We denote 
$\partial \alpha=\{p,q\}$.



\section{2-angled disks}
\LABEL{s:2-angledDisks}

In this section, we investigate a disk called a 2-angled disk
whose boundary contains exactly two white vertices
and consists of edges of label $m$.

Let $\Gamma$ be a chart, $m$ a label of $\Gamma$, $D$ a disk with $\partial D\subset \Gamma_m$, 
and $k$ a positive integer.
If $\partial D$ contains exactly
$k$ white vertices, 
then $D$ is called 
{\it a $k$-angled disk of $\Gamma_m$}. 
Note that 
the boundary $\partial D$ may contain crossings.

Let $\Gamma$ be a chart, and
$m$ a label of $\Gamma$.
An edge of label $m$ is called a {\it feeler} of a $k$-angled disk $D$ of $\Gamma_m$
if the edge intersects $N-\partial D$
where $N$ is a regular neighborhood of $\partial D$ in $D$.

Let $X$ be a set in a chart $\Gamma$.
Let
 $$w(X)=\text{the number of white vertices in $X$.}$$

Let $\Gamma$ be a chart. 
Suppose that an object consists of 
some edges of $\Gamma$, arcs in edges of $\Gamma$ and arcs around white vertices.
Then the object is called a {\it pseudo chart}.

\begin{lemma}
\LABEL{Theorem2AngledDisk}
{\rm (\cite[Corollary 5.8]{ChartAppII})}
Let $\Gamma$ be a minimal chart.
Let $D$ be a $2$-angled disk of $\Gamma_m$ with at most one feeler.
If $w(\Gamma\cap{\rm Int}D)=0$,
then a regular neighborhood of $D$ contains one of two pseudo charts as shown in Fig.~\ref{Fig03}.
\end{lemma}

\begin{figure}
\centerline{\includegraphics{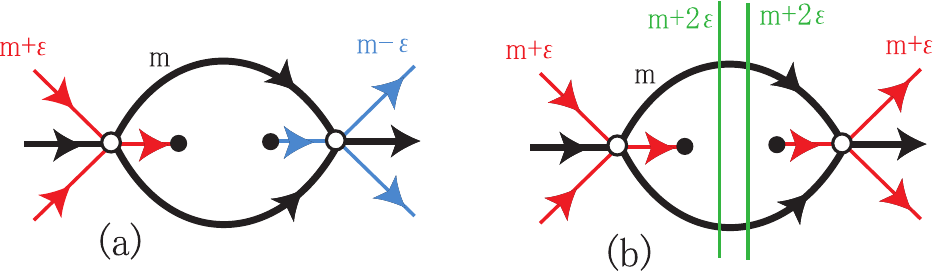}}
\caption{\LABEL{Fig03}
$m$ is a label,
and $\varepsilon\in\{+1,-1\}$.}
\end{figure}

By Lemma~\ref{Theorem2AngledDisk},
we have the following lemma:

\begin{lemma}
\LABEL{2AngledDiskWithOneWhite}
Let $\Gamma$ be a minimal chart,
and $m$ a label of $\Gamma$.
Let $D$ be a $2$-angled disk $D$ of $\Gamma_m$ with 
at most one feeler.
If $D$ satisfies one of the following two conditions,
then $w(\Gamma\cap{\rm Int}D)\ge1$.

\begin{enumerate}
\item[{\rm (a)}]
$D$ has one feeler.
\item[{\rm (b)}]
$D$ has no feeler, and the boundary $\partial D$ is oriented clockwise or anticlockwise.
\end{enumerate}
\end{lemma}

 In our argument  we often construct a chart $\Gamma$. 
On the construction of a chart $\Gamma$, for a white vertex $w\in\Gamma_m$ for some label $m$,  
among the three edges of $\Gamma_m$ 
containing $w$, 
if one of the three edges is a terminal edge 
(see Fig.~\ref{Fig04}(a) and (b)), 
then we remove the terminal edge and
put a black dot at the center of the white vertex  as shown in Fig.~\ref{Fig04}(c).
Namely
Fig.~\ref{Fig04}(c) means 
Fig.~\ref{Fig04}(a) or 
Fig.~\ref{Fig04}(b).
We call the vertex in Fig.~\ref{Fig04}(c) 
a {\it BW-vertex} with respect to $\Gamma_m$.

\begin{figure}[htb]
\centerline{\includegraphics{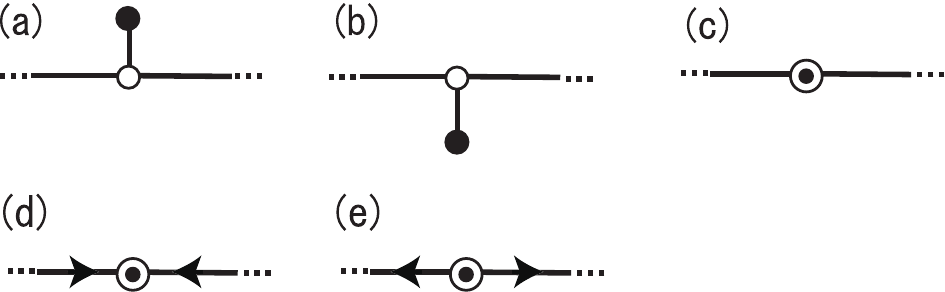}}
\caption{\LABEL{Fig04}
(a),(b) White vertices in terminal edges.
(c),(d),(e) BW-vertices.}
\end{figure}

\begin{lemma}
\LABEL{OriBWvertex}
{\rm (\cite[Lemma 3.1]{ChartAppV})}
In a minimal chart $\Gamma$,
for each BW-vertex in $\Gamma_m$,
the two edges of label $m$ containing the BW-vertex
are oriented inward or outward at the BW-vertex
simultaneously
if each of the two edges is not a terminal edge
$($see Fig.~\ref{Fig04}$($d$)$ and $($e$))$.
\end{lemma}

The three graphs in Fig.~\ref{Fig05}
are examples of graphs in $\Gamma_m$ for a chart $\Gamma$
and a label $m$.
We call 
a {\it $\theta$-curve},
an {\it oval},
a {\it skew $\theta$-curve} the three graphs as shown
in Fig.~\ref{Fig05}(a),(b),(c)
respectively.

\begin{figure}
\centerline{\includegraphics{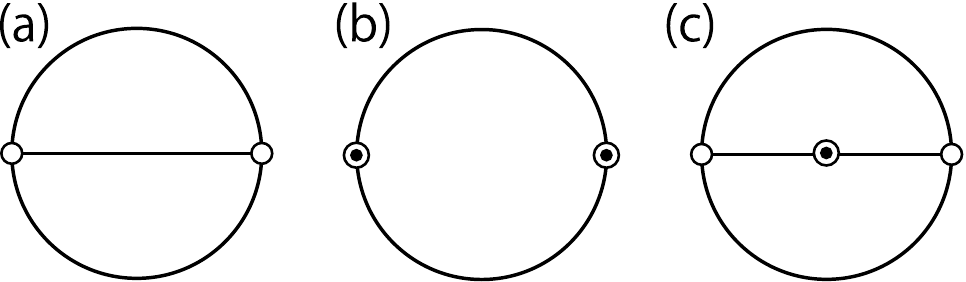}}
\caption{\LABEL{Fig05}
(a) A $\theta$-curve. (b) An oval. 
(c) A skew $\theta$-curve.}
\end{figure}

Let $\Gamma$ be a chart, 
and $m$ a label of $\Gamma$.
Let $L$ be the closure of a connected component 
of the set obtained by taking out 
all the white vertices from $\Gamma_m$.
If $L$ contains at least one white vertex
but does not contain any black vertex,
then $L$ is called an {\it internal edge of label $m$}.
Note that an internal edge may contain a crossing of $\Gamma$.

\begin{lemma}
\LABEL{ThetaCurve}
Let $\Gamma$ be a minimal chart,
and $m$ a label of $\Gamma$.
Let $G$ be a $\theta$-curve in $\Gamma_m$.
Then there exist two connected components
of $S^2-G$ each of which contains at least 
one white vertex.
\end{lemma}

\begin{Proof}
Let $w_1,w_2$ be the white vertices in $G$,
and $e_1,e_2,e_3$ the three internal edges of label $m$ in $G$ such that
$e_1$ is middle at $w_1$.

Let $F_1,F_2$ be the connected components of 
$S^2-G$ with $\partial F_1=e_1\cup e_2$
and $\partial F_2=e_1\cup e_3$
 (see Fig.~\ref{Fig06}).
Then both of $Cl(F_1)$ and $Cl(F_2)$ 
are 2-angled disks without feelers.

Without loss of generality
we can assume that 
$e_1$ is oriented from $w_1$ to $w_2$.
Since $e_1$ is middle at $w_1$,
the two edges $e_2,e_3$ are oriented from $w_2$ to $w_1$
(see Fig.~\ref{Fig06}).
Thus $\partial F_1$ and $\partial F_2$
are oriented clockwise or anticlockwise.
Hence by Lemma~\ref{2AngledDiskWithOneWhite}(b), we have $w(\Gamma\cap F_1)\ge1$
and $w(\Gamma\cap F_2)\ge1$. 
\end{Proof}

\begin{figure}
\centerline{\includegraphics{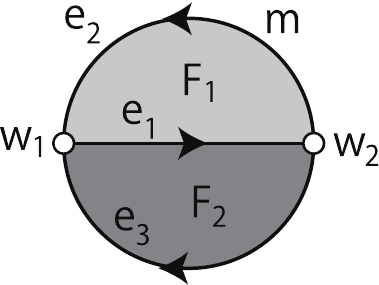}}
\caption{\LABEL{Fig06}
The light gray region is $F_1$,
and the dark gray region is $F_2$.}
\end{figure}

\begin{lemma}
\LABEL{ChartAppLemma67}
{\rm (\cite[Lemma 6.7]{ChartAppII})}
Let $\Gamma$ be a minimal chart, and $m$ a label of $\Gamma$.
Let $D$ be a $2$-angled disk of $\Gamma_m$
with two feelers such that 
$\partial D$ is contained in an oval of label $m$.
Then $w(\Gamma\cap{\rm Int}D)\ge2$.
\end{lemma}

\begin{lemma}
\LABEL{2AngledDiskWithoutFeelers}
Let $\Gamma$ be a minimal chart,
and $m$ a label of $\Gamma$.
Let $D$ be a $2$-angled disk of $\Gamma_m$
without feelers, and
$w_1,w_2$ the white vertices in $\partial D$.
Let $e_1,e_2$ be the internal edges
$($possibly terminal edges$)$ of label $m$
at $w_1,w_2$, respectively,
such that $e_1\not\subset D$
and $e_2\not\subset D$.
Suppose that the two edges $e_1,e_2$ are oriented inward
$($resp. outward$)$ at $w_1,w_2$, 
respectively
$($see Fig.~\ref{Fig07}$($a$)$ and $($b$))$.
Then we have the following.
\begin{enumerate}
\item[{\rm (a)}]
$w(\Gamma\cap{\rm Int}D)\ge1$.
\item[{\rm (b)}]
If $D$ contains an oval of label $m$,
then $w(\Gamma\cap{\rm Int}D)\ge3$.
\end{enumerate}
\end{lemma}

\begin{Proof}
Suppose that $e_1,e_2$ are oriented inward
at $w_1,w_2$, respectively.

We shall show Statement (a).
Let $e_3,e_4$ be the internal edges in $\partial D$.
Without loss of generality
we can assume that $e_3$ is oriented from $w_1$ to $w_2$.
Since both of $e_2$ and $e_3$ are oriented 
inward at $w_2$,
the edge $e_4$ is oriented from $w_2$ to $w_1$ (see Fig.~\ref{Fig07}(c)).
Hence by Lemma~\ref{2AngledDiskWithOneWhite}(b),
we have $w(\Gamma\cap{\rm Int}D)\ge1$.

We shall show Statement (b).
Suppose that
the disk $D$ contains an oval of label $m$,
say $G$.
Let $E$ be the 2-angled disk of $\Gamma_m$
in $D$ with $\partial E\subset G$.
There are three cases:
(i) $E$ has two feelers,
(ii) $E$ has one feeler,
(iii) $E$ has no feeler.

{\bf Case (i).}
Since $\partial E$ is contained in the oval $G$,
by Lemma~\ref{ChartAppLemma67}
we have  $w(\Gamma\cap{\rm Int}E)\ge2$. 
Thus the condition $E\subset {\rm Int}D$
implies $w(\Gamma\cap{\rm Int}D)\ge4\ge3$.

{\bf Case (ii).}
By Lemma~\ref{2AngledDiskWithOneWhite}(a),
we have  $w(\Gamma\cap{\rm Int}E)\ge1$.
Hence we have $w(\Gamma\cap{\rm Int}D)\ge3$.

{\bf Case (iii).}
We use the orientations as shown in 
Fig.~\ref{Fig07}(c).
Let $e_1',e_2'$ be internal edges
(possibly terminal edges)
of label $m\pm1$ at $w_1,w_2$,
respectively,
in $D$.
Since the edges $e_1,e_2$ are oriented inward
at $w_1,w_2$, respectively,

\begin{enumerate}
\item[$(1)$] the edges $e_1',e_2'$ are oriented outward 
at $w_1,w_2$, respectively.
\end{enumerate}
Moreover
since $e_4,e_3$ are oriented inward at $w_1,w_2$, respectively,
neither $e_1'$ nor $e_2'$ is middle at $w_1$ or $w_2$.
Thus by Assumption~\ref{AssumeTerminal},
\begin{enumerate}
\item[$(2)$] neither $e_1'$ nor $e_2'$ 
is a terminal edge.
\end{enumerate}

Let $w_3,w_4$ be the white vertices in the oval $G$.
Since the disk $D$ contains the oval $G$,
by Lemma~\ref{OriBWvertex}
 the two internal edges in $\partial E$
are oriented inward at $w_3$
or oriented outward at $w_3$.
Without loss of generality
we can assume that the two internal edges 
in $\partial E$
are oriented inward at $w_3$.
Since $E$ has no feeler,
the disk $D$ contains the pseudo chart
as shown in Fig.~\ref{Fig07}(d).

If $w(\Gamma\cap({\rm Int}D-E))=0$,
then by (1) and (2)
both of $e_1'$ and $e_2'$ contain $w_4$.
Hence the condition $w(\Gamma\cap({\rm Int}D-E))=0$ implies that
there exist two terminal edges of label $m\pm1$ at $w_3$.
However this contradicts Assumption~\ref{AssumeTerminal}.
Thus $w(\Gamma\cap({\rm Int}D-E))\ge1$.
Therefore $w(\Gamma\cap{\rm Int}D)\ge3$.

Similarly,
we can show Statement (a) and (b) for the case that
 $e_1,e_2$ are oriented outward
at $w_1,w_2$, respectively.
Thus we complete the proof of 
Lemma~\ref{2AngledDiskWithoutFeelers}.
\end{Proof}

\begin{figure}
\centerline{\includegraphics{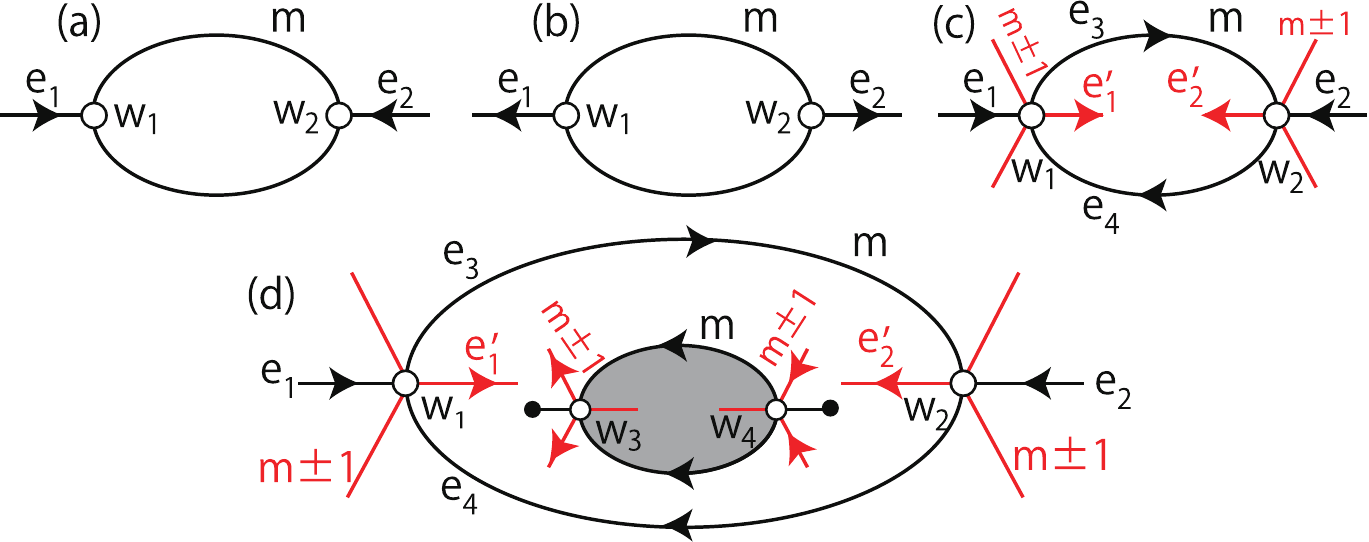}}
\caption{\LABEL{Fig07}
The gray region is the disk $E$. }
\end{figure}



\section{3-angled disks}
\LABEL{s:3-angledDisks}

In this section, we review a lemma about a 3-angled disk.

Let $\Gamma$ and $\Gamma^\prime $ be C-move equivalent charts. 
Suppose that a pseudo chart $X$ of $\Gamma$ is also a pseudo chart of $\Gamma^\prime$. 
Then we say that 
$\Gamma$ is modified to $\Gamma^\prime$ by {\it C-moves keeping $X$ fixed}.
In Fig.~\ref{Fig08},
we give examples of C-moves keeping pseudo charts  fixed.

\begin{figure}[htb]
\centerline{\includegraphics{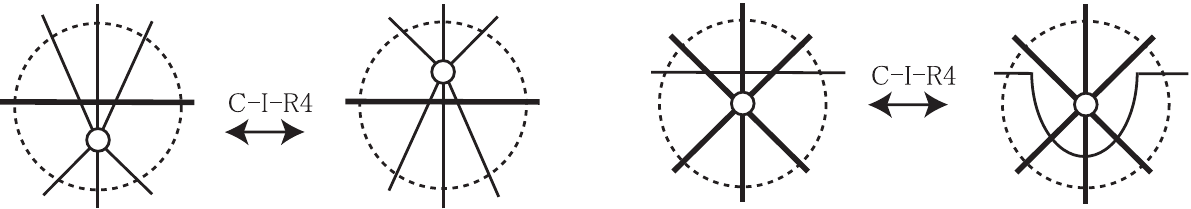}}
\caption{\LABEL{Fig08} 
C-moves keeping thicken figures fixed.}
\end{figure}

Let $X$ be a set in a chart $\Gamma$.
Let
 $$c(X)=\text{the number of crossings in $X$.}$$

Let $D$ be a $k$-angled disk of $\Gamma_m$ for a minimal chart $\Gamma$.
The pair of integers 
$(w(\Gamma\cap{\rm Int}D),c(\Gamma\cap\partial D))$
is called the {\it local complexity 
with respect to $D$},
denoted by $\ell c(D;\Gamma)$.
Let
${\Bbb S}$ be the set of all minimal charts each of which can be moved from $\Gamma$ by C-moves in a regular neighborhood of $D$ keeping $\partial D$ fixed.
The chart $\Gamma$ is said to be 
{\it locally minimal
with respect to $D$}
if its local complexity
with respect to $D$
is minimal
among the charts in ${\Bbb S}$ with respect to 
the lexicographic order.

Let $\Gamma$ be a chart,
and $D$ a $k$-angled disk of $\Gamma_m$.
If any feeler of $D$ of label $m$ is a terminal edge,
then $D$ is called a {\it special} $k$-angled disk.

Let $\Gamma$ be a chart, 
$D$ a disk, and 
$G$ a pseudo chart with $G \subset D$.
Let $r:D\to D$ be a reflection of $D$, and $G^*$ the pseudo chart obtained from $G$ by changing the orientations of all of the edges.
Then the set $\{G,G^*, r(G), r(G^*)\}$ 
is called the {\it RO-family of the pseudo chart $G$}.

\begin{lemma}
{\rm (\cite[Theorem 1.2]{ChartAppIII})}
\LABEL{Theorem3AngledDisk}
Let $\Gamma$ be a minimal chart.
Let $D$ be a special $3$-angled disk of $\Gamma_m$
such that $\Gamma$ is locally minimal
with respect to $D$.
If $w(\Gamma\cap {\rm Int}D)=0$,
then a regular neighborhood of $D$ contains one of the RO-families of the two pseudo charts as shown in Fig.~\ref{Fig09}.
\end{lemma}

\begin{figure}[htb]
\centerline{\includegraphics{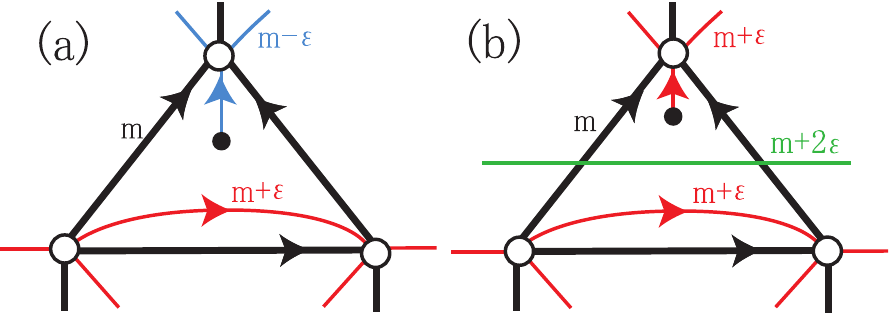}}
\caption{\LABEL{Fig09} 
The 3-angled disks without feelers.
$m$ is a label, $\varepsilon\in\{+1,-1\}$.}
\end{figure}

By Lemma~\ref{Theorem3AngledDisk},
we have the following lemma:

\begin{lemma}
\LABEL{3AngledDiskWithOneWhite}
Let $\Gamma$ be a minimal chart,
and $m$ a label of $\Gamma$.
Let $D$ be a special $3$-angled disk $D$ of $\Gamma_m$ such that
$\Gamma$ is locally minimal chart with respect to $D$.
If $D$ satisfies one of the following two conditions,
then $w(\Gamma\cap{\rm Int}D)\ge1$.

\begin{enumerate}
\item[{\rm (a)}]
$D$ has at least one feeler.
\item[{\rm (b)}]
$D$ has no feeler, and the boundary $\partial D$ is oriented clockwise or anticlockwise.
\end{enumerate}
\end{lemma}



\section{Types of $\Gamma_m$ for a minimal chart of 
type $(m;7)$}
\LABEL{s:TypeGammaM}

In this section, we investigate the graph $\Gamma_m$
for a minimal chart $\Gamma$ of type $(m;7)$.

Let $\Gamma$ be a chart,
and $m$ a label of $\Gamma$. 
A {\it loop} is a simple closed curve in $\Gamma_m$ with exactly one white vertex
(possibly with crossings).

\begin{lemma}
\LABEL{LemmaWithTerminal3}
{\rm (\cite[Lemma 3.2]{ChartAppV})}
Let $\Gamma$ be a minimal chart,
and $m$ a label of $\Gamma$.
Let $G$ be a connected component of $\Gamma_m$.
Then we have the following.
\begin{enumerate}
\item[{\rm (a)}] If $1\le w(G)$, then $2\le w(G)$.
\item[{\rm (b)}] If $1\le w(G)\le 3$
and $G$ does not contain any loop, 
then $G$ is one of three graphs as shown 
in Fig.~\ref{Fig05}.
\end{enumerate}
\end{lemma}

Let $\Gamma$ be a chart, 
and $m$ a label of $\Gamma$.
Let $G_1,G_2,\cdots,G_k$ be all of the connected components of $\Gamma_m$
each of which contains at least one white vertex.
For $i=1,2,\cdots,k$,
we put $n_i=w(G_i)$.
Without loss of generality
we can assume $n_1\ge n_2\ge \cdots\ge n_k$.
Then we say that
$\Gamma_m$ is {\it of type $(n_1,n_2,\cdots,n_k)$.}

\begin{lemma}
\LABEL{LemmaGammaMType(7)}
Let $\Gamma$ be a minimal chart of type $(m;7)$.
Then 
$\Gamma_m$ $($resp. $\Gamma_{m+1})$ is 
of type $(7)$, $(5,2)$, $(4,3)$, or $(3,2,2)$.
\end{lemma}

\begin{Proof}
If $\Gamma_m$ is of type $(n_1,n_2,\cdots,n_k)$,
then by Lemma~\ref{LemmaWithTerminal3}(a)
we have $n_i\ge2$ for $i=1,2,\cdots,k$.
Since $\Gamma$ is of type $(m;7)$,
we have $7=w(\Gamma_m)=n_1+n_2+\cdots+n_k$.

If $k\ge4$,
then 
$$7=n_1+n_2+\cdots+n_k\ge n_1+n_2+n_3+n_4\ge2+2+2+2=8.$$
This is a contradiction.
Hence $1\le k\le 3$.

If $k=1$,
then $n_1=7$.
Hence $\Gamma_m$ is of type $(7)$.
If $k=2$,
then $n_1+n_2=7$, $n_1\ge2$ and $n_2\ge2$.
Hence $\Gamma_m$ is of type $(5,2)$ or $(4,3)$.
If $k=3$,
then $n_1+n_2+n_3=7$, $n_1\ge2$, $n_2\ge2$ and $n_3\ge2$.
Hence $\Gamma_m$ is of type $(3,2,2)$.
Therefore  $\Gamma_m$ is of type $(7)$, $(5,2)$, $(4,3)$ or $(3,2,2)$.

Similarly we can show that
 $\Gamma_{m+1}$ is of type $(7)$, $(5,2)$, $(4,3)$ or $(3,2,2)$.
\end{Proof}

\begin{lemma}$(${\rm \cite[Theorem 1.1]{ChartAppIV}}$)$
\LABEL{LemmaNoLoop}
There is no loop in any minimal chart with exactly seven white vertices.
\end{lemma}

\begin{lemma}
\LABEL{NoType43Type322}
Let $\Gamma$ be a minimal chart of type $(m;7)$.
Then $\Gamma_m$ $($resp. $\Gamma_{m+1})$
is neither of type $(4,3)$ nor of type $(3,2,2)$.
\end{lemma}

\begin{Proof}
Suppose that $\Gamma_m$ is either 
of type $(4,3)$ or of type $(3,2,2)$.
Then there exists a connected component $G$ 
of $\Gamma_m$ with $w(G)=3$.
By Lemma~\ref{LemmaWithTerminal3}(b)
and Lemma~\ref{LemmaNoLoop},
the graph $G$ is a skew $\theta$-curve.
Let $e_1$ be the terminal edge in $G$.
Let $D_1,D_2$ be the closures of connected
components of $S^2-G$ such that $D_1\supset e_1$ and $D_2$ is a 2-angled disk.
Without loss of generality, we can assume that
$\Gamma$ is locally minimal with respect to
$D_1$.
Since $D_1$ is a special 3-angled disk with 
one feeler $e_1$,
by Lemma~\ref{3AngledDiskWithOneWhite}(a)
we have 
\begin{enumerate}
\item[(1)] $w(\Gamma\cap{\rm Int}D_1)\ge1$.
\end{enumerate}

{\bf Claim 1.} $w(\Gamma\cap{\rm Int}D_2)\ge1$.

{\it Proof of Claim~$1$.}
First,
we shall show that $\Gamma$ contains 
the graph as shown in 
Fig.~\ref{Fig10}.
Let $w_1$ be the white vertex in $e_1$.
Let $e_2,e_3$ be the internal edges of label $m$ at $w_1$. 
Let $w_2,w_3$ be the white vertices different 
from $w_1$ in $e_2,e_3$,
respectively.
Without loss of generality
we can assume that the terminal edge $e_1$
is oriented inward at $w_1$.
By Assumption~\ref{AssumeTerminal},
both of $e_2$ and $e_3$ are oriented outward at $w_1$.
Thus the edges $e_2,e_3$ 
are oriented inward at $w_2,w_3$, respectively
(see Fig.~\ref{Fig10}).

Since the two white vertices $w_2,w_3$ are
contained in $\partial D_2$,
by Lemma~\ref{2AngledDiskWithoutFeelers}(a)
we have $w(\Gamma\cap{\rm Int}D_2)\ge1$.
Hence Claim~$1$ holds. {\hfill {$\square$}\vspace{1.5em}}

By (1) and Claim~$1$,
each of ${\rm Int}D_1$ and ${\rm Int}D_2$ 
contains at least one white vertex.
Let $w_4,w_5$ be white vertices in ${\rm Int}D_1,{\rm Int}D_2$,
respectively.
Since $\Gamma$ is of type $(m;7)$,
we have $w_4,w_5\in \Gamma_m$.
Let $G_1,G_2$ be the connected components of
$\Gamma_m$ with $w_4\in G_1$ and $w_5\in G_2$.
Then $\Gamma_m$ contains at least three connected components $G,G_1,G_2$
each of which contains a white vertex.
Since $\Gamma_m$ is of type $(4,3)$ or 
of type $(3,2,2)$,
the set $\Gamma_m$ is of type $(3,2,2)$.
Thus the condition $w(G)=3$ implies 
$w(G_1)=2$ and $w(G_2)=2$.

{\bf Claim 2.} $w(\Gamma\cap{\rm Int}D_2)\ge3$.

{\it Proof of Claim $2$.}
By Lemma~\ref{LemmaWithTerminal3}(b) 
and Lemma~\ref{LemmaNoLoop},
there are two cases:
(i) $G_2$ is a $\theta$-curve,
(ii) $G_2$ is an oval.

{\bf Case (i).}
Let $E_1,E_2$ be the closures of the connected components of $D_2-G_2$ which are 2-angled disks.
Then by Lemma~\ref{ThetaCurve},
one of $E_1$ and $E_2$ contains at least one white vertex in its interior.
Thus $w(\Gamma\cap{\rm Int}D_2)\ge3$.

{\bf Case (ii).}
Since the 2-angled disk $D_2$ contains an oval $G_2$ of label $m$ and since
the two edges $e_2,e_3$ are oriented inward at $w_2,w_3$, respectively
(see Fig.~\ref{Fig10}),
we have $w(\Gamma\cap{\rm Int}D_2)\ge3$
by Lemma~\ref{2AngledDiskWithoutFeelers}(b).
Hence Claim~$2$ holds. {\hfill {$\square$}\vspace{1.5em}}

Since $\Gamma$ is of type $(m;7)$,
we have $w(\Gamma)=7$.
Thus by Claim~2, we have

$$7=w(\Gamma)\ge w(G)+w(G_1)+w(\Gamma\cap{\rm Int}D_2)\ge 3+2+3=8.$$
This is a contradiction.
Therefore $\Gamma_m$ is neither of type $(4,3)$ nor of type $(3,2,2)$.

Similarly we can show that
 $\Gamma_{m+1}$ is neither of type $(4,3)$ nor of type $(3,2,2)$.
Hence we complete the proof of Lemma~\ref{NoType43Type322}.
\end{Proof}

\begin{figure}
\centerline{\includegraphics{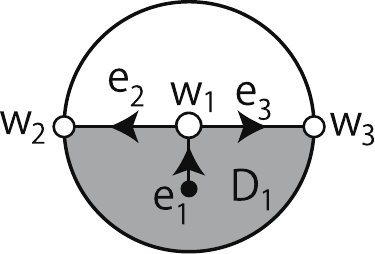}}
\caption{\LABEL{Fig10} 
The gray region is the disk $D_1$.}
\end{figure}



\section{Rings}
\LABEL{s:Rings}

In this section,
we shall show that 
if $\Gamma$ is a minimal chart of type $(m;7)$,
then both of $\Gamma_m$ and $\Gamma_{m+1}$ 
are of type $(5,2)$.

\begin{lemma}
\LABEL{LemmaGammaMType(7)WithRingM-1}
Let $\Gamma$ be a minimal chart of type $(m;7)$.
If $\Gamma$ contains a ring or a non-simple hoop of label $m-1$ $($resp. $m+2)$,
then $\Gamma_m$ $($resp. $\Gamma_{m+1})$ is 
of type $(5,2)$.
\end{lemma}

\begin{Proof}
Suppose that $\Gamma$ contains a ring or a non-simple hoop of label $m-1$.
Let $C$ be the ring or the non-simple hoop.
Let $F_1,F_2$ be the connected components of $S^2-C$.
Then by Assumption~\ref{Ring},
each of $F_1$ and $F_2$ contains at least one white vertex.
Let $w_i$ be the white vertex
in $F_i$ ($i=1,2$).
Since $\Gamma$ is of type $(m;7)$,
we have $w_1,w_2\in \Gamma_m$.
Let $G_1,G_2$ be the connected component of $\Gamma_m$ with $w_i\in G_i$ ($i=1,2$).
Since the curve $C$ of label $m-1$ does not
intersect $\Gamma_m$,
the condition $w_i\in F_i$ ($i=1,2$)
implies $G_i\subset F_i$.
Thus the condition $F_1\cap F_2=\emptyset$
implies $G_1\cap G_2=\emptyset$.
Hence $\Gamma_m$ contains at least two connected components with white vertices.
Thus by Lemma~\ref{LemmaGammaMType(7)} and 
Lemma~\ref{NoType43Type322},
the set $\Gamma_{m}$ is of type $(5,2)$.

Similarly we can show that
if $\Gamma$ contains a ring or a non-simple hoop of label $m+2$,
then the set $\Gamma_{m+1}$ is of 
type $(5,2)$.
\end{Proof}

\begin{lemma}
\LABEL{LemmaType(n)WithRing}
{\rm (\cite[Lemma 6.7]{ChartAppIII})}
Let $\Gamma$ be a minimal chart of type $(m;n_1)$.
Then there exists a ring or a non-simple hoop of label $m-1$ or $m+2$.
\end{lemma}

\begin{lemma}
\LABEL{GammaMType52}
Let $\Gamma$ be a minimal chart of type $(m;7)$.
Then both of $\Gamma_m$ and $\Gamma_{m+1}$
are of type $(5,2)$.
Moreover we have the following.
\begin{enumerate}
\item[{\rm (a)}]
Let $G$ be the connected component of $\Gamma_m$ $($resp. $\Gamma_{m+1})$ with $w(G)=5$.
Then there exists a connected component $F$ of $S^2-G$ with $w(\Gamma\cap F)=2$.
\item[{\rm (b)}]
The other connected components of $S^2-G$
do not contain any white vertices.
\item[{\rm (c)}]
The domain $F$ contains an oval of label $m$
$($resp. $m+1)$
as shown in Fig.~\ref{Fig11}.
\item[{\rm (d)}] $w(\Gamma\cap F')=0$ or $w(\Gamma\cap F')=2$ for
each connected component $F'$ of $S^2-G$.
\end{enumerate}
\end{lemma}

\begin{figure}
\centerline{\includegraphics{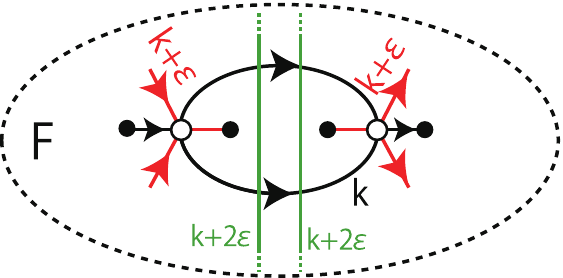}}
\caption{\LABEL{Fig11}
($k=m,\varepsilon=+1$) or ($k=m+1,\varepsilon=-1$). }
\end{figure}

\begin{Proof}
Since $\Gamma$ is of type $(m;7)$,
by Lemma~\ref{LemmaType(n)WithRing}
there exists a ring or a non-simple hoop of label $m-1$ or $m+2$.

Now suppose that there exists a ring or a non-simple hoop of label $m-1$.
Then by Lemma~\ref{LemmaGammaMType(7)WithRingM-1}
the set $\Gamma_m$ is of type $(5,2)$.
Let $G,G'$ be the connected components of $\Gamma_m$ with 
$w(G)=5$ and $w(G')=2$.
Then $G\cap G'=\emptyset$ implies
that the graph $G'$ is contained in
a connected component $F$ of $S^2-G$.
Thus $w(\Gamma\cap F)\ge2$.
If $w(\Gamma\cap F)>2$, then
$$7=w(\Gamma)\ge w(G)+w(\Gamma\cap F)>5+2=7.$$
This is a contradiction.
Hence $w(\Gamma\cap F)=2$.
Thus Statement (a) holds.

Similarly we can show Statement (b).
Statement (d) follows from Statement (a) and (b).

We shall show Statement (c).
By Lemma~\ref{LemmaWithTerminal3}(b) and Lemma~\ref{LemmaNoLoop}
there are two cases:
(i) $G'$ is a $\theta$-curve,
(ii) $G'$ is an oval.

{\bf Case (i).}
By the similar way to the proof of Case (i) of Claim 2 in 
Lemma~\ref{NoType43Type322},
we have $w(\Gamma\cap F)\ge 3$.
This contradicts Statement (d).
Hence Case (i) does not occur.

{\bf Case (ii).}
Let $D$ be the 2-angled disk of $\Gamma_m$
with $\partial D\subset G'$ and $D\subset F$.
Thus the condition $w(\Gamma\cap F)=2$ implies $w(\Gamma\cap{\rm Int}D)=0$.
Hence by Lemma~\ref{Theorem2AngledDisk},
a regular neighborhood of $D$ contains one of 
two pseudo charts as shown in 
Fig.~\ref{Fig03}.
Since $\Gamma$ is of type $(m;7)$,
the two white vertices in $\partial D$
are contained in $\Gamma_m\cap\Gamma_{m+1}$.
Thus a regular neighborhood of $D$ contains 
the pseudo chart as shown in 
Fig.~\ref{Fig03}(b) where $\varepsilon=+1$.
Hence the domain $F$ contains the pseudo chart as shown in Fig.~\ref{Fig11}
 where $k=m$ and $\varepsilon=+1$.
Thus Statement (c) holds.

Moreover the disk $D$ intersects two proper arcs of label $m+2$.
Hence the chart $\Gamma$ contains a ring of label $m+2$.
Thus by Lemma~\ref{LemmaGammaMType(7)WithRingM-1},
the set $\Gamma_{m+1}$ is of type $(5,2)$.
Similarly for the set $\Gamma_{m+1}$
we can show the four statements in this lemma.

Similarly
if there exists a ring or a non-simple hoop of label $m+2$,
we can show this lemma.
Hence we complete the proof of Lemma~\ref{GammaMType52}.
\end{Proof}


\section{Connected components of $\Gamma_m$ with exactly five white vertices}
\LABEL{s:FiveWhiteVertices}

In this section,
we investigate
a connected component of $\Gamma_m$
with exactly five white vertices for a minimal chart $\Gamma$.

\begin{lemma}
\LABEL{LemmaWithTerminal}
{\rm (\cite[Lemma 3.4]{ChartAppVI})}
Let $\Gamma$ be a minimal chart,
and $m$ a label of $\Gamma$.
Let $G$ be a connected component of $\Gamma_m$.
If $w(G)=5$ and $G$ has no loop,
then $G$ is one of nine graphs as shown 
in  Fig.~\ref{Fig12}.
\end{lemma}

\begin{figure}[htb]
\centerline{\includegraphics{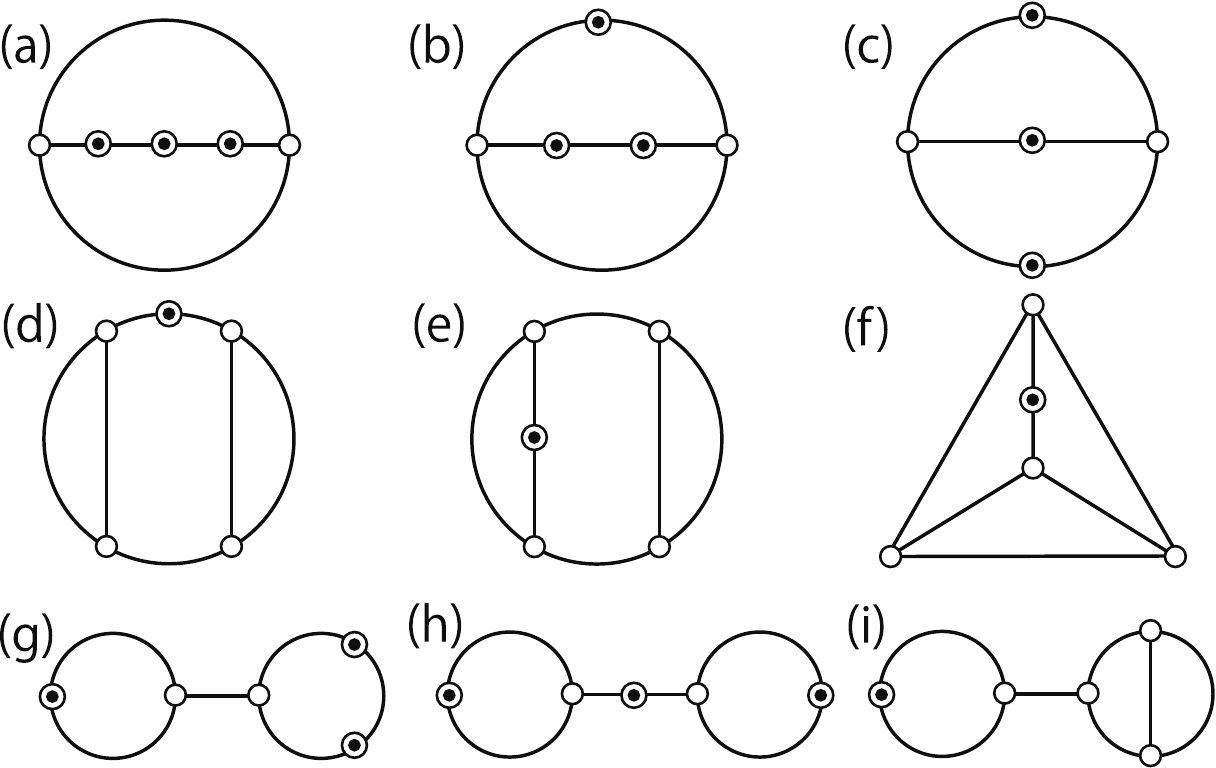}}
\caption{\LABEL{Fig12}
(a),(b),(c) Graphs with three black vertices.
(d),(e),(f) Graphs with one black vertex.
(g),(h) Graphs with three black vertices.
(i) A graph with one black vertex.}
\end{figure}

\begin{lemma}
\LABEL{OriGammaM5}
Let $\Gamma$ be a minimal chart, and $m$ a label of $\Gamma$.
Let $G$ be a connected component of $\Gamma_m$
with $w(G)=5$.
Then we have the following:

\begin{enumerate}
\item[{\rm (a)}]
If $G$ is the graph as shown in Fig.~\ref{Fig12}$($a$)$
$($resp. Fig.~\ref{Fig12}$($b$))$,
then $G$ is one of the RO-family of the graph as shown
in Fig.~\ref{Fig13}$($a$)$
$($resp. Fig.~\ref{Fig13}$($b$))$.
\item[{\rm (b)}]
Suppose that 
$G$ is the graph as shown in Fig.~\ref{Fig12}$($c$)$.
If necessary we move the point at infinity $\infty$,
then $G$ is one of the RO-family of the graph as shown
in Fig.~\ref{Fig13}$($c$)$.
\item[{\rm (c)}]
If $G$ is the graph as shown in Fig.~\ref{Fig12}$($d$)$
$($resp. Fig.~\ref{Fig12}$($e$),($g$),($h$))$,
then $G$ is one of the RO-family of the graph as shown
in Fig.~\ref{Fig13}$($d$)$
$($resp. Fig.~\ref{Fig13}$($e$),($f$),($g$))$.
\end{enumerate}
\end{lemma}

\begin{figure}[htb]
\centerline{\includegraphics{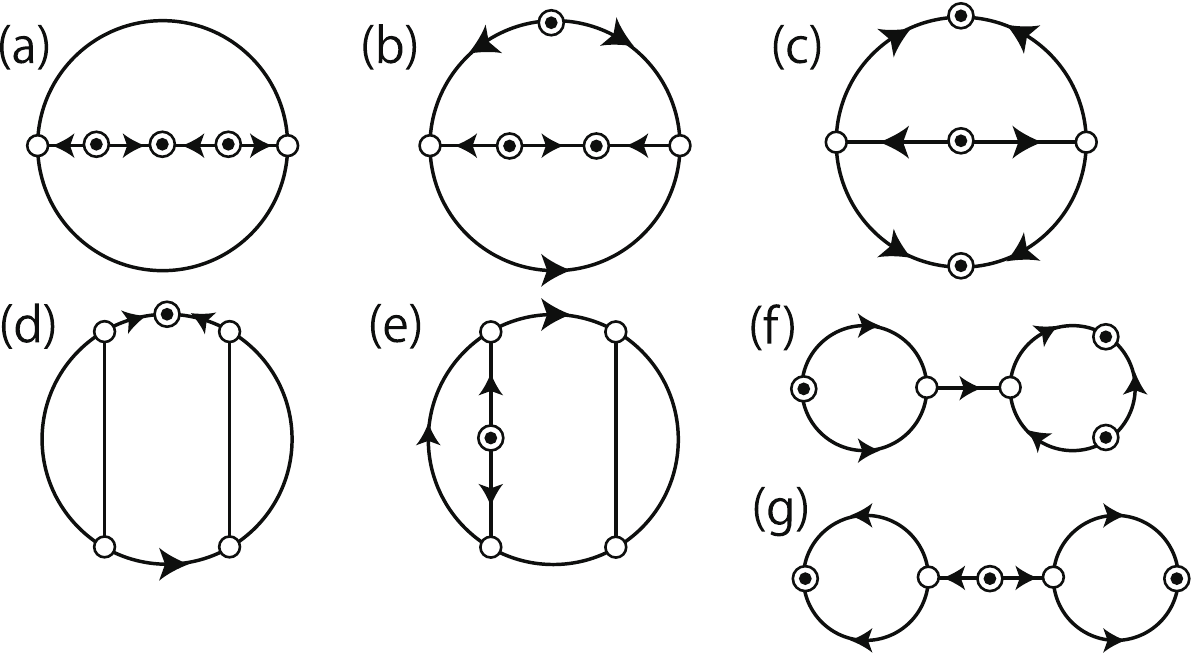}}
\caption{\LABEL{Fig13}
Connected components of $\Gamma_m$ with five white vertices.}
\end{figure}


\begin{Proof}
We show only this lemma for the case that the graph $G$
is one of the three graphs as shown in 
Fig.~\ref{Fig12}(b),(c),(g).

Suppose that $G$ is the graph as shown in Fig.~\ref{Fig12}(b).
We use the notations as shown in Fig.~\ref{Fig14}(a),
where $w_1,w_2,w_3$ are BW-vertices,
$e_1,e_2,\cdots,e_6$ are the six internal edges of label $m$
with $e_1\cap e_2\ni w_1$,
$\partial e_3=\{w_4,w_5\}$,
$\partial e_4=\{w_2,w_4\}$,
$\partial e_5=\{w_2,w_3\}$,
$\partial e_6=\{w_3,w_5\}$.

Without loss of generality we can assume that 
the terminal edge of label $m$ at $w_1$ 
is oriented inward at $w_1$.
Then by Assumption~\ref{AssumeTerminal},
both of $e_1$ and $e_2$ are oriented outward at $w_1$.
Thus
\begin{enumerate}
\item[(1)] $e_2$ is oriented inward at $w_5$.
\end{enumerate}

If necessary we reflect the chart $\Gamma$,
we can assume that 
the edge $e_3$ is oriented from $w_4$ to $w_5$.
Hence by (1),
both of $e_2$ and $e_3$ are oriented inward at $w_5$.
Thus the edge $e_6$ is oriented from $w_5$ to $w_3$.
Hence by Lemma~\ref{OriBWvertex},
the edge $e_5$ is oriented from $w_2$ to $w_3$
and 
the edge $e_4$ is oriented from $w_2$ to $w_4$.
Thus the graph $G$ is the graph as shown in
Fig.~\ref{Fig13}(b).

Suppose that $G$ is the graph as shown in Fig.~\ref{Fig12}(c).
We use the notations as shown in Fig.~\ref{Fig14}(b),
where $w_3,w_4,w_5$ are BW-vertices,
$e_1,e_2,\cdots,e_6$ are the six internal edges of label $m$
with $w_1\in e_1\cap e_2\cap e_3$,
$w_2\in e_4\cap e_5\cap e_6$,
$w_3\in e_1\cap e_4$,
$w_4\in e_2\cap e_5$,
$w_5\in e_3\cap e_6$.

If necessary we move the point at infinity $\infty$
by Assumption~\ref{Infinity},
we can assume that 
\begin{enumerate}
\item[(2)] $e_1$ is middle at $w_1$.
\end{enumerate}
Moreover 
we can assume that 
the edge $e_1$ is oriented inward at $w_1$.
Then by (2),
both of $e_2$ and $e_3$ are oriented outward at $w_1$.
Thus $e_1$ is oriented outward at $w_3$,
and $e_2,e_3$ are oriented inward at $w_4,w_5$,
respectively.
Hence by Lemma~\ref{OriBWvertex},
the edge $e_4$ is oriented outward at $w_3$,
and $e_5,e_6$ are oriented inward at $w_4,w_5$,
respectively.
Thus the graph $G$ is the graph as shown in 
Fig.~\ref{Fig13}(c).

Suppose that $G$ is the graph as shown in Fig.~\ref{Fig12}(g).
We use the notations as shown in Fig.~\ref{Fig14}(c),
where $w_1,w_4,w_5$ are BW-vertices,
$e_1,e_2,\cdots,e_6$ are the six internal edges
of label $m$ with $e_1\cap e_2\ni w_1,w_2$,
$\partial e_3=\{w_2,w_3\}$,
$\partial e_4=\{w_3,w_4\}$,
$\partial e_5=\{w_4,w_5\}$,
$\partial e_6=\{w_3,w_5\}$.

Without loss of generality we can assume that
the terminal edge of label $m$ at $w_1$
is oriented inward at $w_1$.
Then by Assumption~\ref{AssumeTerminal},
both of $e_1,e_2$ are oriented from $w_1$ to $w_2$.
Thus the edge $e_3$ is oriented from $w_2$ to $w_3$.

For the edge $e_4$,
there are two cases:
(i) $e_4$ is oriented from $w_3$ to $w_4$,
(ii) $e_4$ is oriented from $w_4$ to $w_3$.

{\bf Case (i).}
By Lemma~\ref{OriBWvertex},
the edge $e_5$ is  oriented from $w_5$ to $w_4$
and 
the edge $e_6$ is  oriented from $w_5$ to $w_3$.
Thus $G$ is the graph as shown in Fig.~\ref{Fig13}(f).

{\bf Case (ii).}
By Lemma~\ref{OriBWvertex},
the edge $e_5$ is  oriented from $w_4$ to $w_5$
and 
the edge $e_6$ is  oriented from $w_3$ to $w_5$.
Thus by reflecting the chart $\Gamma$,
the graph $G$ is the graph as shown in 
Fig.~\ref{Fig13}(f).
\end{Proof}

\begin{figure}[htb]
\centerline{\includegraphics{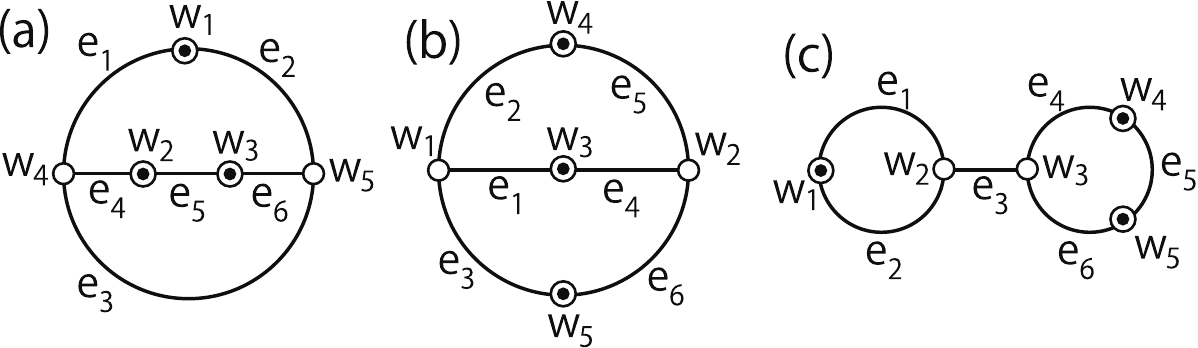}}
\caption{\LABEL{Fig14}
(a) $w_1,w_2,w_3$ are BW-vertices.
(b) $w_3,w_4,w_5$ are BW-vertices.
(c) $w_1,w_4,w_5$ are BW-vertices.
}
\end{figure}



\section{Cases of the graphs as shown in Fig.~\ref{Fig12}(a) and (d)}
\LABEL{s:TypeATypeD}

In this section,
we shall show that
neither $\Gamma_{m}$ nor $\Gamma_{m+1}$ contains
the graphs as shown in Fig.~\ref{Fig12}(a) and (d)
for any minimal chart $\Gamma$ of type $(m;7)$.

\begin{lemma}
\LABEL{2AngledDiskWithoutFeelersType7}
Let $\Gamma$ be a minimal chart of type $(m;7)$.
Let $D$ be a $2$-angled disk of $\Gamma_m$
without feelers, and
$w_1,w_2$ the white vertices in $\partial D$.
Let $e_1,e_2$ be the internal edges
$($possibly terminal edges$)$ of label $m$
at $w_1,w_2$, respectively,
such that $e_1\not\subset D$
and $e_2\not\subset D$.
Then one of $e_1,e_2$ is oriented inward at $w_1$ or $w_2$,
and the other is oriented outward at $w_1$ or $w_2$.
\end{lemma}

\begin{Proof}
Since $\Gamma$ is of type $(m;7)$,
by Lemma~\ref{GammaMType52}
the set $\Gamma_m$ is of type $(5,2)$.
Thus the boundary $\partial D$
is contained in a connected component $G$ of
$\Gamma_m$ with $w(G)=5$ or $w(G)=2$.

If $w(G)=2$,
then by Lemma~\ref{GammaMType52}(c)
 the graph $G$ is an oval as shown in Fig.~\ref{Fig11}
where $k=m$ and $\varepsilon=+1$.
Hence both of $e_1,e_2$ are terminal edges, and
one of $e_1,e_2$ is oriented inward at $w_1$ or $w_2$,
and the other is oriented inward at $w_1$ or $w_2$.

Now suppose $w(G)=5$.
If $e_1,e_2$ are oriented inward 
at $w_1,w_2$, respectively,
then by Lemma~\ref{2AngledDiskWithoutFeelers}(a) we have $w(\Gamma \cap{\rm Int}D)\ge1$.
Hence by Lemma~\ref{GammaMType52}, 
the disk $D$ contains an oval of label $m$
as shown in Fig.~\ref{Fig11}
where $k=m$ and $\varepsilon=+1$.
Thus by Lemma~\ref{2AngledDiskWithoutFeelers}(b), we have $w(\Gamma \cap{\rm Int}D)\ge3$.
This contradicts Lemma~\ref{GammaMType52}(d).
Therefore one of $e_1,e_2$ is oriented outward at $w_1$ or $w_2$.

Similarly
if $e_1,e_2$ are oriented outward at $w_1,w_2$, respectively,
then we have the same contradiction.
Therefore one of  $e_1,e_2$ is oriented inward at $w_1$ or $w_2$.
\end{Proof}


\begin{lemma}
\LABEL{NoTypeA}
Let $\Gamma$ be a minimal chart of type $(m;7)$.
Then neither $\Gamma_m$ nor $\Gamma_{m+1}$
contains the graph as shown in Fig.~\ref{Fig12}$($a$)$.
\end{lemma}

\begin{Proof}
Suppose that $\Gamma_m$ contains the graph
as shown in Fig.~\ref{Fig12}(a), say $G$.
By Lemma~\ref{OriGammaM5}(a),
the graph $G$ is one of the RO-family of the graph
as shown in Fig.~\ref{Fig13}(a).
Without loss of generality
we can assume that $G$ is the graph as shown in 
Fig.~\ref{Fig13}(a).

Let $D$ be the 2-angled disk of $\Gamma_m$
without feelers with $\partial D\subset G$,
and $w_1,w_2$ the two white vertices in $\partial D$.
Let $e_1,e_2$ be the internal edges of label $m$
at $w_1,w_2$, respectively,
with $e_1\not\subset \partial D$ and
$e_2\not\subset \partial D$. 
Then the two edges $e_1,e_2$ are oriented inward at $w_1,w_2$,
respectively (see Fig.~\ref{Fig15}(a)).
However this contradicts
Lemma~\ref{2AngledDiskWithoutFeelersType7}.
Therefore $\Gamma_m$ does not contain the graph as shown in Fig.~\ref{Fig12}(a).

Similarly we can show that  $\Gamma_{m+1}$ does not contain the graph as shown in Fig.~\ref{Fig12}(a).
\end{Proof}

\begin{figure}
\centerline{\includegraphics{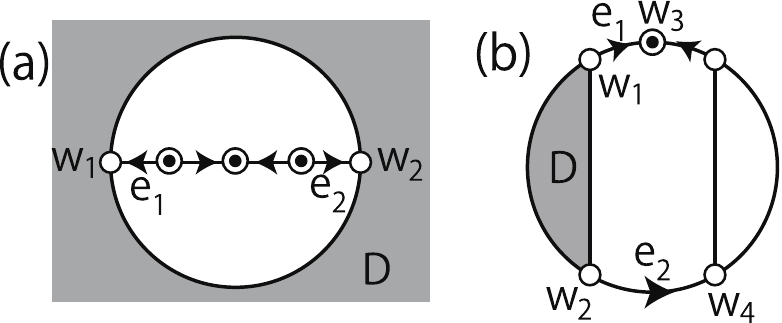}}
\caption{\LABEL{Fig15} 
(a) The graph as shown in Fig.~\ref{Fig12}(a).
(b) The graph as shown in Fig.~\ref{Fig12}(d).
}
\end{figure}


\begin{lemma}
\LABEL{NoTypeD}
Let $\Gamma$ be a minimal chart of type $(m;7)$.
Then neither $\Gamma_m$ nor $\Gamma_{m+1}$
contains the graph as shown in Fig.~\ref{Fig12}$($d$)$.
\end{lemma}

\begin{Proof}
Suppose that $\Gamma_m$ contains the graph
as shown in Fig.~\ref{Fig12}(d), say $G$.
By Lemma~\ref{OriGammaM5}(c),
the graph $G$ is one of the RO-family of the graph
as shown in Fig.~\ref{Fig13}(d).
Without loss of generality
we can assume that $G$ is the graph as shown in 
Fig.~\ref{Fig13}(d).
We use the notations as shown in
Fig~\ref{Fig15}(b),
where $e_1,e_2$ are internal edges of label $m$
with $\partial e_1=\{w_1,w_3\}$ and
$\partial e_2=\{w_2,w_4\}$.

Let $D$ be the 2-angled disks of $\Gamma_m$ without feelers 
with $\partial D\ni w_1,w_2$.
Then the two edges $e_1,e_2$ are oriented outward at $w_1,w_2$,
respectively
 (see Fig.~\ref{Fig15}(b)).
However this contradicts
Lemma~\ref{2AngledDiskWithoutFeelersType7}.
Therefore $\Gamma_m$ does not contain the graph as shown in Fig.~\ref{Fig12}(d).

Similarly we can show that  $\Gamma_{m+1}$ does not contain the graph as shown in Fig.~\ref{Fig12}(d).
\end{Proof}




\section{IO-Calculation}

\LABEL{s:TypeF}

In this section,
we shall show that
neither $\Gamma_{m}$ nor $\Gamma_{m+1}$ contains
the graph as shown in Fig.~\ref{Fig12}(f)
for any minimal chart $\Gamma$ of type $(m;7)$.

Let $\Gamma$ be a chart,
 and $v$ a vertex. 
Let $\alpha$ be a short arc of $\Gamma$ in a small neighborhood of $v$ such that $v$ is an endpoint of $\alpha$. 
If the arc $\alpha$ is oriented to $v$, then $\alpha$ is called {\it an inward arc}, 
and otherwise $\alpha$ is called {\it an outward arc}.

Let $\Gamma$ be an $n$-chart. 
Let $F$ be a closed domain with $\partial F\subset \Gamma_{k-1}\cup\Gamma_{k}\cup \Gamma_{k+1}$ for some label $k$ of $\Gamma$, where $\Gamma_0=\emptyset$ and $\Gamma_{n}=\emptyset$. 
By Condition (iii) for charts,
in a small neighborhood of each white vertex, there are three inward arcs and three outward arcs.
Also in a small neighborhood of each black vertex, there exists only one inward arc or one outward arc.
We often use the following fact, 
when we fix (inward or outward) arcs 
near white vertices and black vertices: 
\begin{enumerate}
\item[$(*)$]
{\it The number of inward arcs contained in $F\cap \Gamma_k$ is equal to the number of outward arcs in $F\cap \Gamma_k$.
}
\end{enumerate}
When we use this fact, 
we say that we use {\it IO-Calculation with respect to $\Gamma_k$ in $F$}.
For example, in a minimal chart $\Gamma$, 
consider the pseudo chart as shown in Fig.~\ref{Fig16} 
where
\begin{enumerate}
\item[(1)]  $v_1,v_2,\cdots,v_5$ are five white vertices in $\Gamma_{k-1}\cap\Gamma_{k}$,
\item[(2)] $D_1$ is a $3$-angled disk of $\Gamma_{k-1}$ without feelers with $\partial D_1\ni v_1,v_2,v_3$,
\item[(3)] $D_1$ contains an oval $G$ of label $k-1$ with $G\ni v_4,v_5$, 
\item[(4)] $e_1,e_2,e_3$ are internal edges 
(possibly terminal edges) of label $k$ 
at $v_1,v_2,v_3$, respectively.
\end{enumerate}
Let $D_2$ be the 2-angled disk of $\Gamma_{k-1}$ in $D_1$ with $\partial D_2\subset G$.
Set $F=Cl(D_1-D_2)$.
If none of $e_1,e_2,e_3$ are middle 
at $v_1,v_2,v_3$, respectively,
then we can show that $w(\Gamma\cap{\rm Int}F)\ge1$.
Suppose $w(\Gamma\cap{\rm Int}F)=0$.
By Assumption~\ref{AssumeTerminal},
\begin{enumerate}
\item[(5)]  none of edges  $e_1,e_2,e_3$ 
are terminal edges. 
\end{enumerate}
Let $e_4',e_4'',e_5',e_5''$ be internal edges
 (possibly terminal edges) of label $k$ 
such that $e_4',e_4''$ are oriented inward at $v_4$ and 
$e_5',e_5''$ are oriented outward at $v_5$.
Since none of edges  $e_4',e_4'',e_5',e_5''$
are middle at $v_4$ or $v_5$,
by Assumption~\ref{AssumeTerminal}
\begin{enumerate}
\item[(6)] none of edges  $e_4',e_4'',e_5',e_5''$ are terminal edges. 
\end{enumerate}
If the edges $e_1,e_2,e_3$ are oriented 
inward at $v_1,v_2,v_3$, respectively,
then 
by (5) and (6) 
the number of inward arcs in $F\cap \Gamma_{k}$ is five,  
but the number of outward arcs in $F\cap \Gamma_{k}$ is two. 
This contradicts the fact $(*)$. 
Similarly for the other cases
 we have the same contradiction.
Thus $w(\Gamma\cap{\rm Int}F)\ge1$.
Instead of the above argument, 
we just say that if none of $e_1,e_2,e_3$ are middle at $v_1,v_2,v_3$, respectively, then
\begin{enumerate}
\item[]
{\it we have $w(\Gamma\cap{\rm Int}F)\ge1$ 
by IO-Calculation with respect to $\Gamma_{k}$ in $F$.}
\end{enumerate}

\begin{figure}
\centerline{\includegraphics{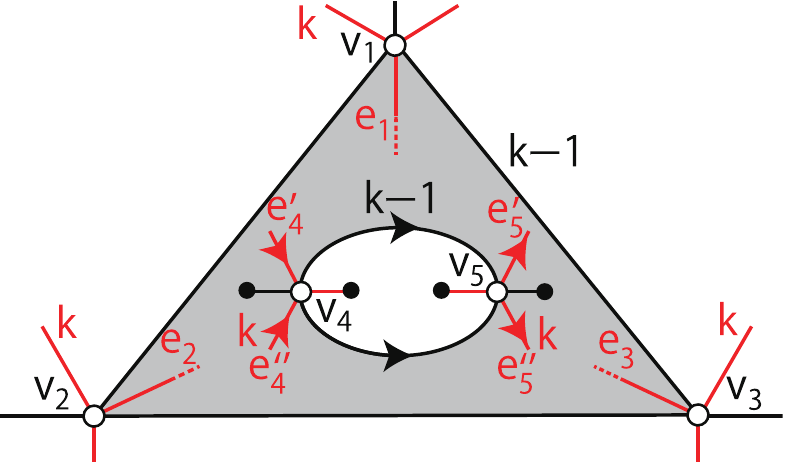}}
\caption{\LABEL{Fig16} The gray region is $F$.}
\end{figure}

\begin{lemma}
\LABEL{No3angledDiskClockwise}
Let $\Gamma$ be a minimal chart of type $(m;7)$.
Then there does not exist a $3$-angled disk  $D$ of $\Gamma_m$ without feelers such that
$\partial D$ is oriented clockwise or anticlockwise.
\end{lemma}

\begin{Proof}
Suppose that there exists a $3$-angled disk $D$ of $\Gamma_m$ without feelers such that
$\partial D$ is oriented clockwise or anticlockwise.
Without loss of generality, we can assume 
that $\Gamma$ is locally minimal with respect
to $D$.
Then by Lemma~\ref{3AngledDiskWithOneWhite}(b) we have 
 $w(\Gamma\cap{\rm Int}D)\ge1$.

Let $G$ be the connected component of $\Gamma_m$ with $G\supset \partial D$.
Since $\Gamma$ is of type $(m;7)$,
by Lemma~\ref{GammaMType52}
the graph $\Gamma_m$ is of type $(5,2)$.
Thus we have $w(G)=5$.
Hence by Lemma~\ref{GammaMType52},
the condition $w(\Gamma\cap{\rm Int}D)\ge1$
implies that the disk $D$ contains an oval of label $m$
 as shown in Fig.~\ref{Fig11}
where $k=m$ and $\varepsilon=+1$.

Let $w_1,w_2,w_3$ be the white vertices in
$\partial D$,
and $e_1,e_2,e_3$ internal edges (possibly
terminal edges) of label $m+1$ at $w_1,w_2,w_3$, respectively, in $D$.
Let $D'$ be the 2-angled disk of $\Gamma_m$ in $D$ without feelers.
Since $\partial D$ is oriented clockwise or
anticlockwise,
none of $e_1,e_2,e_3$ are middle at
$w_1,w_2,w_3$, respectively.
Thus considering as $D_1=D$, $D_2=D'$, $F=Cl(D-D')$ and $k=m+1$ in the example of IO-Calculation in this section,
we have $w(\Gamma\cap({\rm Int}D-D'))\ge1$.
Hence $w(\Gamma\cap{\rm Int}D)\ge3$.
This contradicts Lemma~\ref{GammaMType52}(d).
Therefore  there does not exist a $3$-angled disk $D$ of $\Gamma_m$ without feelers such that
$\partial D$ is oriented clockwise or anticlockwise.
\end{Proof}


\begin{lemma}
\LABEL{NoTypeF}
Let $\Gamma$ be a minimal chart of type $(m;7)$.
Then neither $\Gamma_m$ nor $\Gamma_{m+1}$
contains the graph as shown in Fig.~\ref{Fig12}$($f$)$.
\end{lemma}

\begin{Proof}
Suppose that $\Gamma_m$ contains the graph
as shown in Fig.~\ref{Fig12}(f).
We use the notations as shown in 
Fig.~\ref{Fig17}(a),
where $w_1$ is the BW-vertex,
and $e_1',e_1'',e_2',e_2'',e_3',e_3'',e_4$ are 
seven internal edges of
label $m$ with
$e_1'\cap e_1''\ni w_1$,
$e_1'\cap e_2'\cap e_2''\ni w_2$,
$e_1''\cap e_3'\cap e_3''\ni w_3$, and
$\partial e_4=\{w_4,w_5\}$.
Let $D_1,D_2$ be the 3-angled disks of $\Gamma_m$ without feelers 
with $\partial D_1=e_3'\cup e_3''\cup e_4$ 
and 
$\partial D_2=e_2'\cup e_2''\cup e_4$.

Without loss of generality
we can assume that
the terminal edge of label $m$ at $w_1$
is oriented inward at $w_1$. 
Then by Assumption~\ref{AssumeTerminal},
both of $e_1',e_1''$ are oriented outward at $w_1$ 
(see Fig.~\ref{Fig17}(a)). 
For the edge $e_1'$,
there are two cases:
(i) $e_1'$ is middle at $w_2$,
(ii) $e_1'$ is not middle at $w_2$.

{\bf Case (i).}
By Condition (iii) of the definition of a chart,
both of $e_2',e_2''$ are oriented outward at $w_2$.
If necessary we reflect the chart $\Gamma$,
we can assume that
\begin{enumerate}
\item[(1)]
the edge $e_4$ is oriented from $w_4$ to $w_5$.
\end{enumerate}
Since both of $e_2'',e_4$ are oriented inward at $w_5$,
\begin{enumerate}
\item[(2)]
the edge $e_3''$ is oriented from $w_5$ to $w_3$.
\end{enumerate}
Moreover,
since both of $e_1'',e_3''$ are oriented inward at $w_3$,
the edge $e_3'$ is oriented from $w_3$ to $w_4$ 
(see Fig.~\ref{Fig17}(b)).
Hence by (1) and (2),
the boundary $\partial D_1$ is oriented clockwise or anticlockwise.
This contradicts Lemma~\ref{No3angledDiskClockwise}.
Hence Case (i) does not occur.

{\bf Case (ii).}
One of $e_2',e_2''$ is oriented inward at $w_2$, and the other is oriented outward at $w_2$.
If necessary we reflect the chart $\Gamma$,
we can assume that 
the edge $e_2'$ is oriented inward at $w_2$, and the edge $e_2''$ is oriented outward at $w_2$.
Looking at the 3-angled disk $D_2$,
by Lemma~\ref{No3angledDiskClockwise}
the edge $e_4$ is oriented from $w_4$ to $w_5$.
By the same way to Case (i),
the edge $e_3''$ is oriented from $w_5$ to $w_3$, and the edge $e_3'$ is oriented from $w_3$ to $w_4$ (see Fig.~\ref{Fig17}(c)).
Thus the boundary $\partial D_1$ is oriented clockwise or anticlockwise.
This contradicts Lemma~\ref{No3angledDiskClockwise}.
Hence Case (ii) does not occur.

Therefore the both cases do not occur.
Hence  $\Gamma_{m}$ does not contain the graph as shown in Fig.~\ref{Fig12}(f).

Similarly we can show that  $\Gamma_{m+1}$ does not contain the graph as shown in Fig.~\ref{Fig12}(f).
\end{Proof}

\begin{figure}
\centerline{\includegraphics{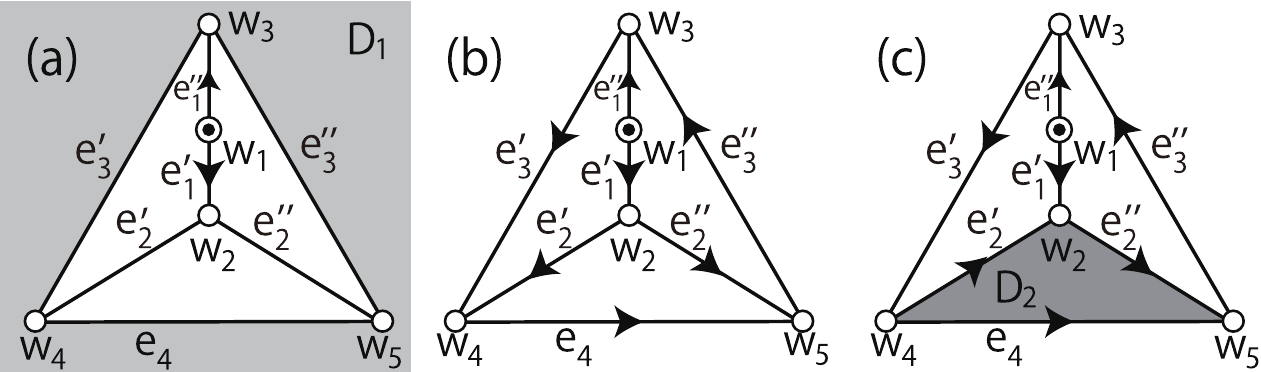}}
\caption{\LABEL{Fig17} 
The graphs as shown in Fig.~\ref{Fig12}(f).
The light gray region is the disk $D_1$.
The dark gray region is the disk $D_2$.}
\end{figure}


\section{Lenses}
\LABEL{s:TypeI}

In this section,
we shall show that
neither $\Gamma_{m}$ nor $\Gamma_{m+1}$ contains
the graph as shown in Fig.~\ref{Fig12}(i)
for any minimal chart $\Gamma$ of type $(m;7)$.

Let $\Gamma$ be a chart. 
Let $D$ be a disk 
such that 
\begin{enumerate}
\item[(1)] the boundary $\partial D$ consists of an internal edge $e_1$ of label $m$ and an internal edge $e_2$ of label ${m+1}$, and 
\item[(2)] any edge containing a white vertex in $e_1$ does not intersect the open disk Int$D$.
\end{enumerate}
Note that $\partial D$ may contain crossings.
Let $w_1$ and $w_2$ be the white vertices in $e_1$. 
If the disk $D$ satisfies one of the following conditions, then $D$ is called  {\it a lens of type $(m,m+1)$}
(see Fig.~\ref{Fig18}):
\begin{enumerate}
	\item[(i)] Neither $e_1$ nor $e_2$ contains a middle arc. 
	\item[(ii)] One of the two edges $e_1$ and $e_2$ contains middle arcs at both white vertices $w_1$ and $w_2$ simultaneously.
\end{enumerate}

\begin{figure}[htb]
\centerline{\includegraphics{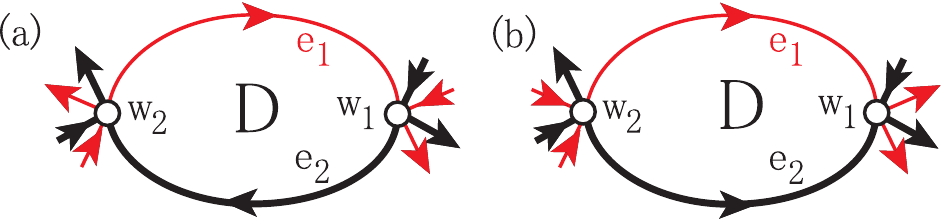}}
\caption{\LABEL{Fig18}
Lenses.}
\end{figure}

\begin{lemma}{\rm (\cite[Corollary 1.3]{ChartAppII})}
\LABEL{NoLens}
 There is no lens in any minimal chart with 
at most seven white vertices.
\end{lemma}

\begin{lemma}
\LABEL{No2angledDiskClockwise}
Let $\Gamma$ be a minimal chart of type $(m;7)$.
Then there does not exist a $2$-angled disk  $D$ of $\Gamma_m$ without feelers such that
$\partial D$ is oriented clockwise or anticlockwise.
\end{lemma}

\begin{Proof}
Suppose that there exists a $2$-angled disk $D$ of $\Gamma_m$ without feelers such that
$\partial D$ is oriented clockwise or anticlockwise.

Let $w_1,w_2$ be the white vertices 
in $\partial D$, and 
$e_1,e_2$ internal edges 
(possibly terminal edges) of label $m$ 
at $w_1,w_2$, respectively,
with $e_1\not\subset D$ and 
$e_2\not\subset D$.
By Lemma~\ref{2AngledDiskWithoutFeelersType7},
one of $e_1,e_2$ is oriented inward at $w_1$ or $w_2$, and 
the other is oriented outward at $w_1$ or $w_2$.
Without loss of generality we can assume that
$e_1$ is oriented inward at $w_1$, and 
$e_2$ is oriented outward at $w_2$.
Let $e_1',e_2'$ be internal edges 
(possibly terminal edges) of label $m+1$
at $w_1,w_2$, respectively, in $D$.
Then 
\begin{enumerate}
\item[(1)] 
$e_1'$ is oriented outward at $w_1$, and 
$e_2'$ is oriented inward at $w_2$.
\end{enumerate}

Let $G$ be the connected component of $\Gamma_m$ with $G\supset \partial D$.
Since $\Gamma_m$ is of type $(5,2)$ 
by Lemma~\ref{GammaMType52},
we have $w(G)=5$ or $w(G)=2$.

If $w(G)=2$, then by Lemma~\ref{GammaMType52}(c)
the graph $G$ is an oval of label $m$ and 
the two internal edges in $\partial D$
are oriented from $w_1$ to $w_2$
or are oriented from $w_2$ to $w_1$
(see Fig.~\ref{Fig11}).
This contradicts that $\partial D$ is oriented clockwise or anticlockwise.
Thus $w(G)=5$.

Since $\partial D$ is oriented clockwise or
anticlockwise,
by Lemma~\ref{2AngledDiskWithOneWhite}(b)
the interior of $D$ contains at least
one white vertex.
Thus by Lemma~\ref{GammaMType52},
the condition $w(G)=5$ implies that
the disk $D$ contains an oval $G'$ of label $m$ as shown in Fig.~\ref{Fig11} where $k=m$ and $\varepsilon=+1$
(see 
Fig.~\ref{Fig19}(a)).
Moreover by Lemma~\ref{GammaMType52}(d),
we have $w(\Gamma\cap{\rm Int}D)=2$.
Let $D'$ be the 2-angled disk of $\Gamma_m$
in $D$ with $\partial D'\subset G'$.
Then the condition $w(\Gamma\cap{\rm Int}D)=2$
implies  
\begin{enumerate}
\item[(2)] $w(\Gamma\cap({\rm Int}D-D'))=0$.
\end{enumerate}

We use the notations as shown in 
Fig.~\ref{Fig19}(a),
where
$w_3,w_4$ are white vertices in the oval $G'$,
$e_3',e_3''$ are internal edges 
(possibly terminal edges) of label $m+1$
oriented inward at $w_3$, and
$e_4',e_4''$ are internal edges 
(possibly terminal edges) of label $m+1$
oriented outward at $w_4$.
By Assumption~\ref{AssumeTerminal},
\begin{enumerate}
\item[(3)] none of $e_3',e_3'',e_4',e_4''$
are terminal edges.
\end{enumerate}
Since $\partial D$ is oriented clockwise
or anticlockwise,
neither $e_1'$ nor $e_2'$ contains 
a middle arc at $w_1$ or $w_2$.
Thus by Assumption~\ref{AssumeTerminal},
\begin{enumerate}
\item[(4)] neither $e_1'$ nor $e_2'$
is a terminal edge.
\end{enumerate}
Hence by (1),(2),(3) and (4),
we have $e_1'=e_3'$ or $e_1'=e_3''$ or $e_1'=e_2'$. Thus
\begin{enumerate}
\item[(5)] $e_1'=e_3'$, $e_2'=e_4'$, and $e_3''=e_4''$ (see Fig.~\ref{Fig19}(b)), or
\item[(6)] $e_1'=e_3''$, $e_2'=e_4''$, and $e_3'=e_4'$ (see Fig.~\ref{Fig19}(c)), or
\item[(7)] $e_1'=e_2''$, $e_3'=e_4'$, and $e_3''=e_4''$ (see Fig.~\ref{Fig19}(d)).
\end{enumerate}
For each case,
there exists a lens of type $(m,m+1)$.
This contradicts Lemma~\ref{NoLens}.
Therefore there does not exist a $2$-angled disk  $D$ of $\Gamma_m$ without feelers such that
$\partial D$ is oriented clockwise or anticlockwise.
\end{Proof}

\begin{figure}
\centerline{\includegraphics{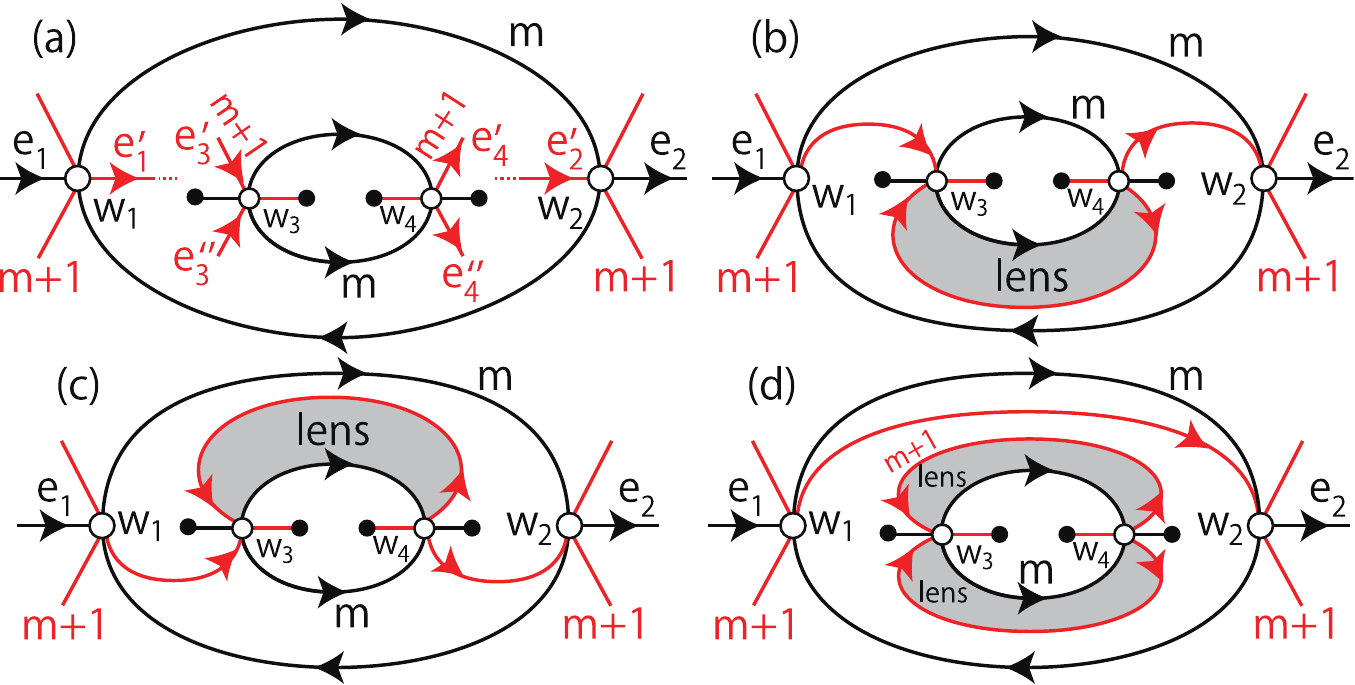}}
\caption{\LABEL{Fig19}
The gray regions are lenses. }
\end{figure}

\begin{lemma}
\LABEL{NoTypeI}
Let $\Gamma$ be a minimal chart of type $(m;7)$.
Then neither $\Gamma_m$ nor $\Gamma_{m+1}$
contains the graph as shown in Fig.~\ref{Fig12}$($i$)$.
\end{lemma}

\begin{Proof}
Suppose that $\Gamma_m$ contains the graph
as shown in Fig.~\ref{Fig12}(i).
We use the notations as shown in 
Fig.~\ref{Fig20}(a),
where $e_1,e_2,e_3,e_4$ are internal edges of
label $m$ with
$\partial e_1=\partial e_2=\{w_1,w_2\}$,
$\partial e_3=\{w_1,w_3\}$,
$\partial e_4=\{w_2,w_3\}$.

Let $D_1,D_2$ be the 2-angled disk
and the 3-angled disk of $\Gamma_m$ without feelers 
with $\partial D_1=e_1\cup e_2$ and 
$\partial D_2=e_1\cup e_3\cup e_4$. 
Without loss of generality
we can assume that
\begin{enumerate}
\item[(1)]
the edge $e_1$ is oriented from $w_1$ to $w_2$.
\end{enumerate}

If the edge $e_2$ is oriented from $w_2$ to $w_1$,
then the boundary $\partial D_1$ is oriented clockwise or 
anticlockwise.
This contradicts Lemma~\ref{No2angledDiskClockwise}. Thus
\begin{enumerate}
\item[(2)]
the edge $e_2$ is oriented from $w_1$ to $w_2$.
\end{enumerate}

By (1) and (2),
the both edges $e_1$ and $e_2$ are oriented outward at $w_1$.
Hence the edge $e_3$ is oriented inward at $w_1$.
Thus the edge $e_3$ is oriented from $w_3$ to $w_1$.
Similarly
by (1) and (2),
the edge $e_4$ is oriented from $w_2$ to $w_3$.
Thus the boundary $\partial D_2$
is oriented clockwise or anticlockwise
(see Fig.~\ref{Fig20}(b)).
This contradicts Lemma~\ref{No3angledDiskClockwise}.
Therefore $\Gamma_m$ does not contain the graph as shown in Fig.~\ref{Fig12}(i).

Similarly we can show that  $\Gamma_{m+1}$ does not contain the graph as shown in Fig.~\ref{Fig12}(i).
\end{Proof}

\begin{figure}
\centerline{\includegraphics{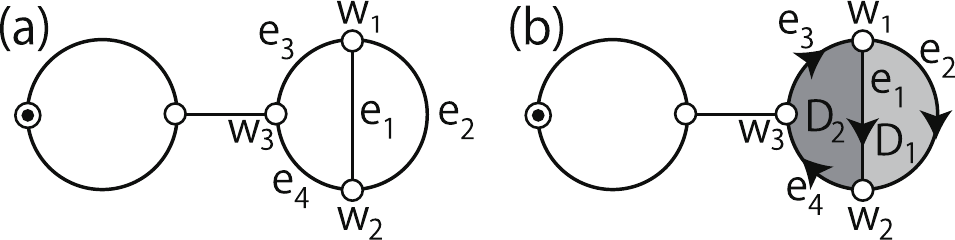}}
\caption{\LABEL{Fig20} 
The graphs as shown in Fig.~\ref{Fig12}(i).
The light gray region is the disk $D_1$.
The dark gray region is the disk $D_2$.
}
\end{figure}


\section{Case of the graph as shown in Fig.~\ref{Fig12}(h)}
\LABEL{s:TypeH}

In this section,
we shall show that
neither $\Gamma_{m}$ nor $\Gamma_{m+1}$ contains
the graph as shown in Fig.~\ref{Fig12}(h)
for any minimal chart $\Gamma$ of type $(m;7)$.

\begin{lemma}
\LABEL{No2angledDiskOneFeeler}
Let $\Gamma$ be a minimal chart of type $(m;7)$.
Then there does not exist a special $2$-angled disk of $\Gamma_m$ with exactly one feeler.
\end{lemma}

\begin{Proof}
Suppose that there exists a special $2$-angled disk $D$ of $\Gamma_m$ with exactly one feeler.
Let $e_1$ be the feeler of $D$,
and $w_1$ the white vertex in $e_1$.
Let  $e_1',e_1''$ be internal edges 
(possibly terminal edges) of label $m+1$
 at $w_1$ in $D$.
Without loss of generality we can assume that
the terminal edge $e_1$ is oriented inward at $w_1$.
Then by Assumption~\ref{AssumeTerminal},
\begin{enumerate}
\item[(1)] 
 $e_1',e_1''$ are oriented inward at $w_1$, 
\item[(2)] neither $e_1'$ nor $e_1''$
is a terminal edge.
\end{enumerate}

Since the 2-angled disk $D$ has exactly one feeler,
 by Lemma~\ref{GammaMType52}
\begin{enumerate}
\item[(3)] the boundary $\partial D$ is contained in a connected component
of $\Gamma_m$ with exactly five white vertices.
\end{enumerate}

Since the 2-angled disk $D$ has one feeler $e_1$,
by Lemma~\ref{2AngledDiskWithOneWhite}(a)
we have $w(\Gamma\cap{\rm Int}D)\ge1$.
Thus by (3) and Lemma~\ref{GammaMType52},
the disk $D$ contains an oval $G$ of label $m$ as shown in Fig.~\ref{Fig11} where $k=m$ and $\varepsilon=+1$
(see 
Fig.~\ref{Fig21}).

Let $D'$ be the 2-angled disk of $\Gamma_m$
in $D$ with $\partial D'\subset G$.
We use the notations as shown in 
Fig.~\ref{Fig21},
where $w_2$ is a white vertex in $G$, and
$e_2',e_2''$ are internal edges 
(possibly terminal edges) of label $m+1$
oriented inward at $w_2$.
By Assumption~\ref{AssumeTerminal},
neither $e_2'$ nor $e_2''$
is a terminal edge.
Hence by (1) and (2),
there are four internal edges 
$e_1',e_1'',e_2',e_2''$ of label $m+1$
oriented inward at 
$w_1,w_1,w_2,w_2$,
respectively, in $D$.
Thus by IO-Calculation with respect to
$\Gamma_{m+1}$ in $Cl(D-D')$,
we have $w(\Gamma\cap({\rm Int}D-D'))\ge1$.
Hence $w(\Gamma\cap{\rm Int}D)\ge3$.
This contradicts Lemma~\ref{GammaMType52}(d).
Therefore there does not exist a special $2$-angled disk 
of $\Gamma_m$ with exactly feeler.
\end{Proof}

\begin{figure}
\centerline{\includegraphics{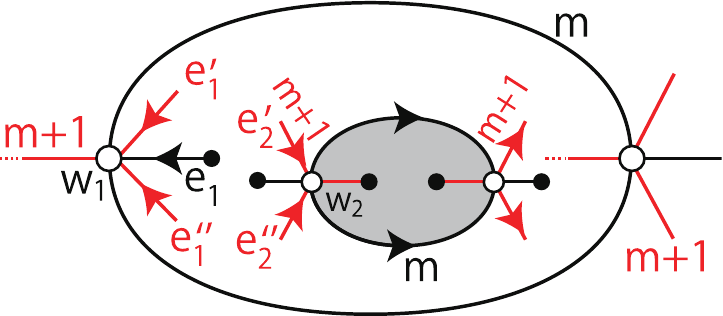}}
\caption{\LABEL{Fig21} 
The gray region is the disk $D'$.}
\end{figure}

\begin{lemma}
\LABEL{NoTypeH}
Let $\Gamma$ be a minimal chart of type $(m;7)$.
Then neither $\Gamma_m$ nor $\Gamma_{m+1}$
contains the graph as shown in Fig.~\ref{Fig12}$($h$)$.
\end{lemma}

\begin{Proof}
Suppose that $\Gamma_m$ contains the graph
as shown in Fig.~\ref{Fig12}(h), say $G$.
By Lemma~\ref{OriGammaM5}(c),
the graph $G$ is one of the RO-family of the graph
as shown in Fig.~\ref{Fig13}(g).
Without loss of generality
we can assume that
the graph $G$ is the graph
as shown in Fig.~\ref{Fig13}(g).

Let $D_1,D_2$ be the 2-angled disks of $\Gamma_m$
with $D_1\cap D_2=\emptyset$,
$\partial D_1\subset G$,
$\partial D_2\subset G$.
Then by Lemma~\ref{No2angledDiskOneFeeler},
both of $D_1$ and $D_2$
have no feelers.
Thus the chart $\Gamma$ contains the
pseudo chart as shown in Fig.~\ref{Fig22}.

Since the graph $G$ is the graph as shown in Fig.~\ref{Fig12}(h),
the set $S^2-G$ consists of three connected component of $S^2-G$.
Two of three are ${\rm Int}D_1$ and 
${\rm Int}D_2$.
Let $F$ be the last connected component 
of $S^2-G$.

We use the notations as shown in Fig.~\ref{Fig22},
where $e_1',e_1'',e_2',e_2'',e_3',e_3''$
are six internal edges 
(possibly terminal edges)
of label $m+1$ oriented inward at 
$w_1,w_1,w_2,w_2,w_3,w_3$,
respectively.
Since none of them are middle at $w_1,w_2$ or $w_3$,
by Assumption~\ref{AssumeTerminal}
none of them
are terminal edges.
Hence by IO-Calculation with respect to $\Gamma_{m+1}$ in $Cl(F)$,
we have $w(\Gamma\cap F)\ge1$.
Thus by Lemma~\ref{GammaMType52},
the domain $F$ contains
an oval $G'$ of label $m$ as shown in
Fig.~\ref{Fig11},
where $k=m$ and $\varepsilon=+1$.

Let $D_3$ be the 2-angled disk of $\Gamma_m$
without feelers whose boundary is 
contained in $G'$.
Again by IO-Calculation with respect to $\Gamma_{m+1}$ in $Cl(F-D_3)$,
we have $w(\Gamma\cap (F-D_3))\ge1$.
Hence $w(\Gamma\cap F)\ge3$.
This contradicts Lemma~\ref{GammaMType52}(d).
Therefore $\Gamma_m$ does not contain the graph as shown in Fig.~\ref{Fig12}(h).

Similarly we can show that $\Gamma_{m+1}$ does not contain the graph as shown in Fig.~\ref{Fig12}(h).
\end{Proof}

\begin{figure}
\centerline{\includegraphics{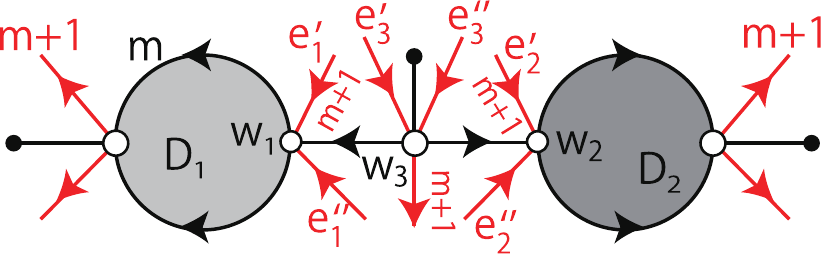}}
\caption{\LABEL{Fig22}
The six edges $e_1',e_1'',e_2',e_2'',e_3',e_3''$ are
oriented inward at $w_1,w_1,w_2,w_2,w_3,w_3$, respectively.}
\end{figure}


\section{Case of the graph as shown in Fig.~\ref{Fig12}(e)}
\LABEL{s:TypeE}

In this section,
we shall show that
neither $\Gamma_{m}$ nor $\Gamma_{m+1}$ contains
the graph as shown in Fig.~\ref{Fig12}(e)
for any minimal chart $\Gamma$ of type $(m;7)$.

\begin{lemma}
\LABEL{Special3angledDiskOneOrTwoFeelers}
Let $\Gamma$ be a minimal chart of type $(m;7)$.
Then any special $3$-angled disk of $\Gamma_m$ has no feelers. 
\end{lemma}

\begin{Proof}
Suppose that there exists a special $3$-angled disk $D$ of $\Gamma_m$ with at least one feeler.
By Lemma~\ref{3AngledDiskWithOneWhite}(a),
we have $w(\Gamma\cap{\rm Int D})\ge1$.
Thus by Lemma~\ref{GammaMType52},
the disk $D$ contains the oval $G$ of label $m$ as shown in Fig.~\ref{Fig11}
where $k=m$ and $\varepsilon=+1$
(see 
Fig.~\ref{Fig23}(a)).
We use the notations as shown in 
Fig.~\ref{Fig23}(a),
where
\begin{enumerate}
\item[(1)] $D'$ is the 2-angled disk of 
$\Gamma_m$ in $D$,
\item[(2)] 
 $e_1',e_1''$ are internal edges 
(possibly terminal edges) of label $m+1$
oriented inward at $w_1$.
\end{enumerate}
By Assumption~\ref{AssumeTerminal},
\begin{enumerate}
\item[(3)] none of $e_1',e_1'',e_2',e_2''$
are terminal edges.
\end{enumerate}

Since $D$ has at least one feeler,
we can assume that the white vertex $w_3$
is contained in the feeler $e_3$ (see 
Fig.~\ref{Fig23}(a)).
Moreover we can assume that
the terminal edge $e_3$ is oriented inward at $w_3$.
Then by Assumption~\ref{AssumeTerminal},
\begin{enumerate}
\item[(4)] both of $e_3',e_3''$ are oriented inward at $w_3$,
\item[(5)] neither $e_3'$ nor $e_3''$ is a terminal edge, and
\item[(6)] the two internal edges $e,e'$ of label $m$
are oriented outward at $w_3$.
\end{enumerate}

If the 3-angled disk $D$ has three feelers,
then there exists a connected component
of $\Gamma_m$ with exactly three white vertices and three terminal edges.
This contradicts Lemma~\ref{LemmaWithTerminal3}(b).
Thus there are two cases:
(i) $D$ has one feeler (see 
Fig.~\ref{Fig23}(b)),
(ii) $D$ has two feelers (see 
Fig.~\ref{Fig23}(c)).

{\bf Case (i).}
If necessary
we reflect the chart $\Gamma$,
we can assume that the edge $e''$ of label $m$ is 
oriented from $w_4$ to $w_5$.
Thus by (6)
the edge $e_5$ is oriented inward at $w_5$
(see 
Fig.~\ref{Fig23}(b)).
Since $e_1',e_1'',e_3',e_3'',e_5$ are
oriented inward at $w_1,w_1,w_3,w_3,w_5$,
respectively by (2) and (4), and
since none of $e_1',e_1'',e_3',e_3''$
are terminal edges by (3) and (5),
we have $w(\Gamma\cap({\rm Int}D-D'))\ge1$
by IO-Calculation with respect to 
$\Gamma_{m+1}$ in $Cl(D-D')$.
Hence  $w(\Gamma\cap{\rm Int}D)\ge3$.
This contradicts Lemma~\ref{GammaMType52}(d).
Thus Case (i) does not occur.

{\bf Case (ii).}
Without loss of generality
we can assume that the white vertex $w_4$
is contained in a feeler $e_4$.
Since the edge $e$ is oriented inward at $w_4$ by (6), 
the edge $e''$ is oriented inward at $w_4$ 
by Assumption~\ref{AssumeTerminal}.
Thus the edge $e''$ is oriented outward at $w_5$ 
(see Fig.~\ref{Fig23}(c)).
Since 
the edge $e'$ is oriented inward at $w_5$ by (6),
the edge $e_5$ is not middle at $w_5$.
Hence by Assumption~\ref{AssumeTerminal},
\begin{enumerate}
\item[(7)] the edge $e_5$ is not a terminal edge. 
\end{enumerate}

Let $e_4',e_4''$ be internal edges
(possibly terminal edges) of label $m+1$
at $w_4$ in $D$.
By Assumption~\ref{AssumeTerminal},
none of $e_4',e_4''$ are terminal edges.
Thus by (3), (5) and (7),
none of the nine edges
$e_1',e_1'',e_2',e_2'',e_3',e_3'',e_4',e_4'',e_5$ 
are terminal edges.
Hence we have $w(\Gamma\cap({\rm Int}D-D'))\ge1$ by IO-Calculation with respect to
$\Gamma_{m+1}$ in $Cl(D-D')$.
Thus $w(\Gamma\cap{\rm Int}D)\ge3$.
This contradicts Lemma~\ref{GammaMType52}(d).
Hence Case (ii) does not occur.

Therefore the both cases do not occur.
Hence any special $3$-angled disk of $\Gamma_m$ has no feelers.
\end{Proof}

\begin{figure}
\centerline{\includegraphics{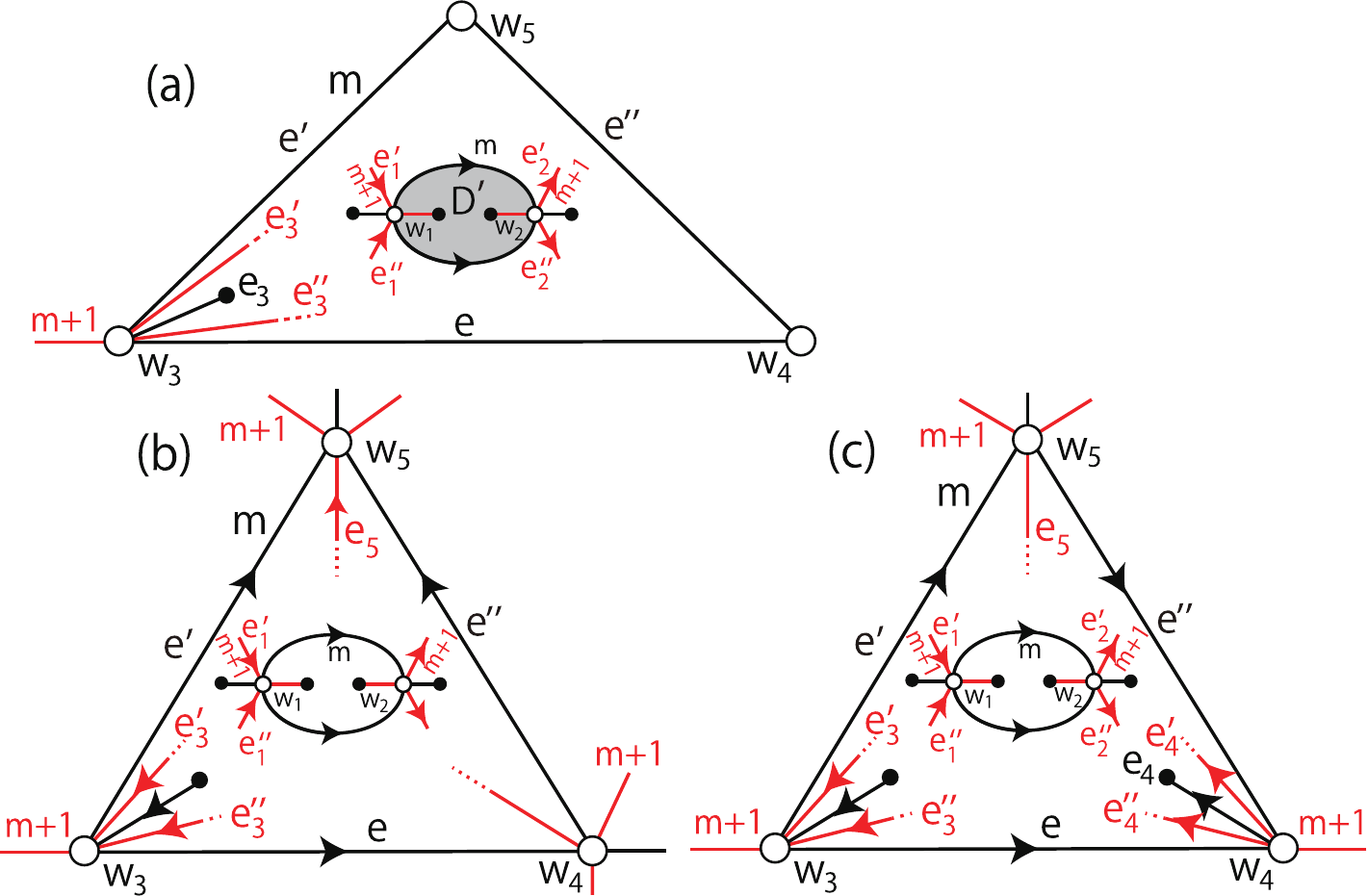}}
\caption{\LABEL{Fig23}
(a) A 3-angled disk. 
(b) A 3-angled disk with one feeler.
(c) A 3-angled disk with two feelers.
}
\end{figure}

\begin{lemma}
\LABEL{NoTypeE}
Let $\Gamma$ be a minimal chart of type $(m;7)$.
Then neither $\Gamma_m$ nor $\Gamma_{m+1}$
contains the graph as shown in Fig.~\ref{Fig12}$($e$)$.
\end{lemma}

\begin{Proof}
Suppose that $\Gamma_m$ contains the graph
as shown in Fig.~\ref{Fig12}(e), say $G$.
By Lemma~\ref{OriGammaM5}(c),
the graph $G$ is one of the RO-family of the graph
as shown in Fig.~\ref{Fig13}(e).
Without loss of generality
we can assume that
the graph $G$ is the graph
as shown in Fig.~\ref{Fig13}(e).

We use the notations as shown in 
Fig.~\ref{Fig24}(a),
where $w_5$ is the BW-vertex, and
$e_1,e_2,\cdots,e_7$ are seven internal edges of
label $m$ with
$\partial e_1=\partial e_2=\{w_1,w_2\}$,
$\partial e_3=\{w_1,w_3\}$,
$\partial e_4=\{w_2,w_4\}$,
$e_5\cap e_6 \ni w_5$,
$\partial e_7=\{w_3,w_4\}$, and
\begin{enumerate}
\item[(1)]
the edge $e_3$ is oriented inward at $w_1$.
\end{enumerate}

{\bf Claim.}
The chart $\Gamma$ contains the pseudo chart as shown in 
Fig.~\ref{Fig24}(b).

{\it Proof of Claim.}
First we show that both of the edges $e_1,e_2$ 
are oriented outward at $w_1$.
Since $e_1\cup e_2$ bounds a 2-angled disk without feelers,
by Lemma~\ref{No2angledDiskClockwise}
either
 both of the edges $e_1,e_2$ 
are oriented inward at $w_1$ or
are oriented outward at $w_1$.
If both of the edges $e_1,e_2$ 
are oriented inward at $w_1$,
then by (1)
there are three edges $e_1,e_2,e_3$ of label $m$
oriented inward at $w_1$.
This contradicts the definition of a chart.
Hence both of the edges $e_1,e_2$ 
are oriented outward at $w_1$.

Since both of the edges $e_1,e_2$ 
are oriented inward at $w_2$,
the edge $e_4$ is oriented from $w_2$ to $w_4$.

Let $D_1$ be the 3-angled disk of $\Gamma_m$
with at most one feeler and
$\partial D_1=e_5\cup e_6\cup e_7$
(see Fig.~\ref{Fig24}(a)).
By Lemma~\ref{Special3angledDiskOneOrTwoFeelers},
the 3-angled disk $D_1$ has no feelers.
Hence the chart $\Gamma$ contains the pseudo chart as shown in Fig.~\ref{Fig24}(b).
Thus Claim holds. {\hfill {$\square$}}

Let $D_2$ be the 5-angled disk of $\Gamma_m$ with one feeler such that $\partial D_2=e_1\cup e_3\cup e_4\cup e_5\cup e_6$.
Let $D_3$ be the 4-angled disk of $\Gamma_m$ without feelers
such that 
$\partial D_3=e_2\cup e_3\cup e_4\cup e_7$.
We use the notations as shown in Fig.~\ref{Fig24}(b),
where 
\begin{enumerate}
\item[(2)]$e_1',e_4',e_5',e_5''$
are internal edges (possibly terminal edges)
of label $m+1$ oriented inward at 
$w_1,w_4,w_5,w_5$ in $D_2$,
respectively, 
\item[(3)]  $e_2',e_3',e_4''$
are internal edges (possibly terminal edges)
of label $m+1$ oriented outward at 
$w_2,w_3,w_4$ in $D_3$,
respectively.
\end{enumerate}
Since none of the six edges $e_1',e_5',e_5'',e_2',e_3',e_4''$ are middle at
$w_1,w_5,w_5,w_2,w_3,w_4$, respectively,
by Assumption~\ref{AssumeTerminal}
none of the six edges are terminal edges.
Thus by applying IO-Calculation with respect to $\Gamma_{m+1}$ in $D_2$ and $D_3$,
we have 
$w(\Gamma\cap{\rm Int}D_2)\ge1$
and $w(\Gamma\cap{\rm Int}D_3)\ge1$
by (2) and (3).
However by Lemma~\ref{GammaMType52}(b),
$w(\Gamma\cap{\rm Int}D_2)=0$ or 
$w(\Gamma\cap{\rm Int}D_3)=0$.
This is a contradiction.
Therefore $\Gamma_m$ does not contain the graph as shown in Fig.~\ref{Fig12}(e).

Similarly we can show that $\Gamma_{m+1}$ does not contain the graph as shown in Fig.~\ref{Fig12}(e).
\end{Proof}

\begin{figure}
\centerline{\includegraphics{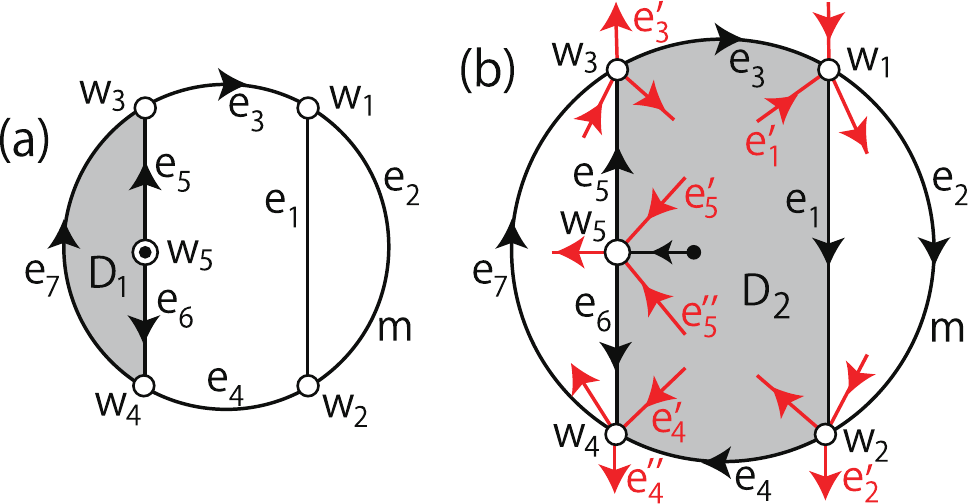}}
\caption{\LABEL{Fig24}
(a) The gray region is the disk $D_1$.
(b) The gray region is the disk $D_2$.}
\end{figure}


\section{Case of the graph as shown in Fig.~\ref{Fig12}(b)}

\LABEL{s:TypeB}

In this section,
we shall show that
neither $\Gamma_{m}$ nor $\Gamma_{m+1}$ contains
the graph as shown in Fig.~\ref{Fig12}(b)
for any minimal chart $\Gamma$ of type $(m;7)$.

\begin{lemma}
\LABEL{Special4Or5angledDiskTwoFeeler}
Let $\Gamma$ be a minimal chart of type $(m;7)$.
Let $D$ be a special $4$-angled 
$($or $5$-angled$)$ disk of $\Gamma_m$ 
with exactly two feelers $e_1,e_2$.
Let $w_1,w_2$ be the white vertices in $e_1,e_2$,
respectively.
Then one of $e_1$ and $e_2$ is oriented inward at $w_1$ or $w_2$ and
the other is oriented outward at 
$w_1$ or $w_2$. 
\end{lemma}

\begin{Proof}
Suppose that both of $e_1,e_2$ are oriented 
inward at $w_1,w_2$,
respectively.
Let $e_1',e_1'',e_2',e_2''$ be internal edges
(possibly terminal edges) of label $m+1$ in $D$ at $w_1,w_1,w_2,w_2$, respectively.
By Assumption~\ref{AssumeTerminal},
the four edges $e_1',e_1'',e_2',e_2''$ are oriented inward at $w_1,w_1,w_2,w_2$, respectively
(see Fig~\ref{Fig25}), and
 none of $e_1',e_1'',e_2',e_2''$ are terminal edges.
Hence by IO-Calculation with respect to $\Gamma_{m+1}$ in $D$,
we have  $w(\Gamma\cap{\rm Int}D)\ge1$.
Thus by Lemma~\ref{GammaMType52},
the disk $D$ contains an oval of label $m$
as shown in Fig.~\ref{Fig11},
where $k=m$ and $\varepsilon=+1$.

Let $D'$ be the 2-angled disk of $\Gamma_m$
in $D$ without feelers.
By IO-Calculation with respect to $\Gamma_{m+1}$ in $Cl(D-D')$,
we have $w(\Gamma\cap({\rm Int}D-D'))\ge1$.
Hence $w(\Gamma\cap{\rm Int}D)\ge3$.
This contradicts Lemma~\ref{GammaMType52}(d).
Therefore one of $e_1,e_2$ is oriented outward at $w_1$ or $w_2$.

Similarly we can show that  one of $e_1,e_2$ is oriented inward at $w_1$ or $w_2$.
Hence we complete the proof of Lemma~\ref{Special4Or5angledDiskTwoFeeler}.
\end{Proof}

\begin{figure}
\centerline{\includegraphics{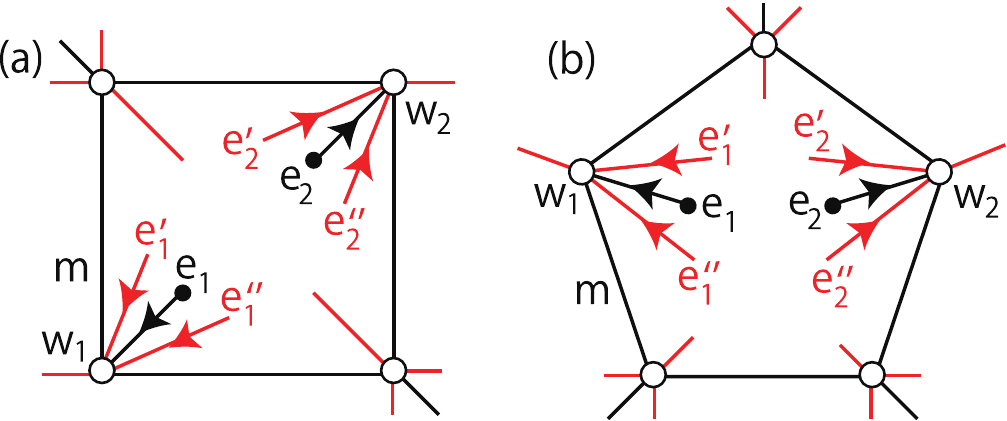}}
\caption{\LABEL{Fig25}
(a) A 4-angled disk with two feelers.
(b) A 5-angled disk with two feelers. }
\end{figure}

\begin{lemma}
\LABEL{Special5angledDiskThreeFeeler}
Let $\Gamma$ be a minimal chart of type $(m;7)$.
Then there does not exist  a special $5$-angled 
 disk of $\Gamma_m$ with exactly three feelers. 
\end{lemma}

\begin{Proof}
Let $D$ be a special 5-angled disk of $\Gamma_m$,
and $w_1,w_2,\cdots,w_5$ the five white vertices 
lying on $\partial D$ in this order.
Let $e_1,e_2,\cdots,e_5$ be the five internal edges of label $m$
in $\partial D$
with $\partial e_i=\{w_i,w_{i+1}\}$ ($i=1,2,3,4$)
and $\partial e_5=\{w_5,w_1\}$.

Suppose that $D$ has exactly three feelers $e,e',e''$.
Since $D$ is special,
all of  $e,e',e''$ are terminal edges.
Without loss of generality we can assume that 
$w_1\in e$, $w_2\in e'$, and
\begin{enumerate}
\item[(1)] one of $w_3,w_4,w_5$ is contained in the feeler $e''$.
\end{enumerate}
Moreover we can assume that the terminal edge $e$
is oriented inward at $w_1$.
Then by Assumption~\ref{AssumeTerminal},
the terminal edge $e$ is middle at $w_1$.
Thus both of $e_1$ and $e_5$ are oriented outward at $w_1$.
Hence $e_1$ is oriented inward at $w_2$.
Thus by Assumption~\ref{AssumeTerminal},
the terminal edge $e'$ is oriented outward at $w_2$, and
\begin{enumerate}
\item[(2)] $e_2$ is oriented inward at $w_2$.
\end{enumerate}
Hence we have the pseudo chart as shown in 
Fig.~\ref{Fig26}(a), where
\begin{enumerate}
\item[(3)] $e_1',e_1''$ are internal edges 
(possibly terminal edges) of label $m+1$ 
oriented inward at $w_1$.
\end{enumerate}
Thus by Assumption~\ref{AssumeTerminal},
\begin{enumerate}
\item[(4)] neither $e_1'$ nor $e_1''$ is a terminal edge.
\end{enumerate}

By (1),
there are three cases:
(i) $w_3\in e''$ (see Fig.~\ref{Fig26}(b)),
(ii) $w_4\in e''$ (see Fig.~\ref{Fig26}(c)),
(iii) $w_5\in e''$.

{\bf Case (i).}
By (2) and Assumption~\ref{AssumeTerminal},
the terminal edge $e''$ is middle at $w_3$ and
oriented inward at $w_3$.
Hence 
\begin{enumerate}
\item[(5)] there exist internal edges $e_3',e_3''$  
(not terminal edges) of label $m+1$ in $D$
oriented inward at $w_3$.
\end{enumerate}
Moreover the edge $e_3$ is oriented from $w_3$ to $w_4$.

If necessary we reflect the chart $\Gamma$,
we can assume that the edge $e_4$ is oriented 
from $w_4$ to $w_5$.
Since $e_5$ is oriented from $w_1$ to $w_5$ 
(see Fig.~\ref{Fig26}(a)),
both of $e_4$ and $e_5$ are oriented inward at $w_5$.
Thus there exists an internal edge $e_5'$
(possibly terminal edge) of label $m+1$ in $D$ oriented inward at $w_5$
(see Fig.~\ref{Fig26}(b)).
Hence by (3) and (5),
the five edges $e_1',e_1'',e_3',e_3'',e_5'$
are oriented inward at $w_1,w_1,w_3,w_3,w_5$,
respectively.
Since none of $e_1',e_1'',e_3',e_3''$ are terminal edges
by (4) and (5),
we have $w(\Gamma\cap{\rm Int}D)\ge1$ by IO-Calculation
with respect to $\Gamma_{m+1}$ in $D$.
Hence by Lemma~\ref{GammaMType52},
the disk $D$ contains an oval of label $m$ as shown in 
Fig.~\ref{Fig11},
where $k=m$ and $\varepsilon=+1$.

Let $D'$ be the 2-angeled disk of $\Gamma_m$ in $D$
without feelers.
By IO-Calculation with respect to $\Gamma_{m+1}$ in $Cl(D-D')$,
we have $w(\Gamma\cap({\rm Int}D-D'))\ge1$.
Thus $w(\Gamma\cap{\rm Int}D)\ge3$.
This contradicts Lemma~\ref{GammaMType52}(d).
Hence Case (i) does not occur.

{\bf Case (ii).}
If necessary we change the orientation of all the edges and
we reflect the chart $\Gamma$,
we can assume that the terminal edge $e''$ is oriented 
inward at $w_4$.
Thus by Assumption~\ref{AssumeTerminal},  
\begin{enumerate}
\item[(6)] there exist internal edges $e_4',e_4''$  
(not terminal edges) of label $m+1$ in $D$
oriented inward at $w_4$, and
\end{enumerate}
both of $e_3,e_4$ are oriented outward at $w_4$.
Since $e_5$ is oriented from $w_1$ to $w_5$
(see Fig.~\ref{Fig26}(a)),
both of $e_4$ and $e_5$ are oriented inward at $w_5$.
Thus there exists an internal edge $e_5'$
(possibly terminal edge) of label $m+1$ in $D$ oriented inward at $w_5$
(see Fig.~\ref{Fig26}(c)).
Hence by (3) and (6),
the five edges $e_1',e_1'',e_4',e_4'',e_5'$
are oriented inward at $w_1,w_1,w_4,w_4,w_5$,
respectively.
Since none of $e_1',e_1'',e_4',e_4''$ are terminal edges
by (4) and (6),
we can show that $w(\Gamma\cap{\rm Int}D)\ge3$ by
the similar way of the proof of Case (i).
This contradicts Lemma~\ref{GammaMType52}(d).
Thus Case (ii) does not occur.

{\bf Case (iii).}
We can show that Case (iii) does not occur by
the similar way of the proof of Case (i).

Therefore all the three cases do not occur.
Hence there does not exist  a special $5$-angled 
 disk of $\Gamma_m$ with exactly three feelers. 
\end{Proof}

\begin{figure}
\centerline{\includegraphics{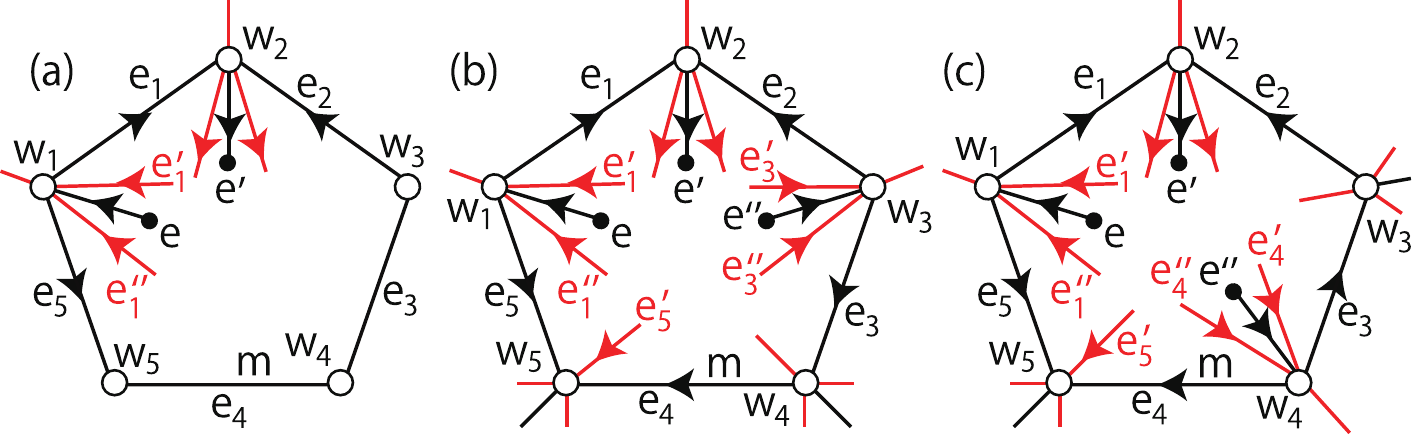}}
\caption{\LABEL{Fig26} 
(a) The white vertices $w_1,w_2$ are contained in feelers $e,e'$.
(b),(c) 5-angled disks with three feelers.}
\end{figure}

\begin{lemma}
\LABEL{NoTypeB}
Let $\Gamma$ be a minimal chart of type $(m;7)$.
Then neither $\Gamma_m$ nor $\Gamma_{m+1}$
contains the graph as shown in Fig.~\ref{Fig12}$($b$)$.
\end{lemma}

\begin{Proof}
Suppose that $\Gamma_m$ contains the graph
as shown in Fig.~\ref{Fig12}(b), say $G$.
By Lemma~\ref{OriGammaM5}(a),
the graph $G$ is one of the RO-family of the graph
as shown in Fig.~\ref{Fig13}(b).
Without loss of generality
we can assume that
the graph $G$ is the graph
as shown in Fig.~\ref{Fig13}(b).

We use the notations as shown in Fig.~\ref{Fig13}(b)
and Fig.~\ref{Fig14}(a),
where $w_1,w_2,w_3$ are BW-vertices,
and $e_2,e_4$ are internal edges of label $m$
with $\partial e_2=\{w_1,w_5\}$ and 
$\partial e_4=\{w_2,w_4\}$.
Let $e,e',e''$ be the terminal edges of label $m$
at $w_1,w_2,w_3$, respectively.
Then
\begin{enumerate}
\item[(1)]  $e,e'$ are oriented inward at $w_1,w_2$, 
respectively.
\end{enumerate}

Now the graph $G$ divides $S^2$ into three disks.
Let $D_1,D_2,D_3$ be the three disks
such that
$D_1$ is a 3-angled disk,
$D_2$ is a 4-angled disk,
$D_3$ is a 5-angled disk.
By Lemma~\ref{Special3angledDiskOneOrTwoFeelers},
the 3-angled disk $D_1$ has no feeler.
Thus the terminal edge $e$ is contained in $D_3$.
Hence looking at the 5-angled disk $D_3$,
by Lemma~\ref{Special5angledDiskThreeFeeler}
we have $e'\not\subset D_3$ or $e''\not\subset D_3$.
Hence $e'\subset D_2$ or $e''\subset D_2$.
Thus there are three cases:
(i) $e'\subset D_2$, $e''\subset D_3$ 
(see Fig.~\ref{Fig27}(a)),
(ii) $e'\subset D_3$, $e''\subset D_2$,
(iii) $e'\cup e''\subset D_2$ (see Fig.~\ref{Fig27}(b)).

{\bf Case (i).}
We use the notations as shown in Fig.~\ref{Fig27}(a),
where
\begin{enumerate}
\item[(2)] $e_2',e_2''$ are internal edges 
(possibly terminal edges) of label $m+1$ 
oriented inward at $w_2$ in $D_2$, 
\item[(3)] $e_3',e_3'',e_5'$ are internal edges 
(possibly terminal edges) of label $m+1$ 
oriented outward at $w_3,w_3,w_5$ in $D_3$,
respectively, and
\end{enumerate}
none of $e_2',e_2'',e_3',e_3'',e_5'$ are middle 
at $w_2,w_2,w_3,w_3,w_5$,
respectively.
Thus by Assumption~\ref{AssumeTerminal},
\begin{enumerate}
\item[(4)] none of the five edges $e_2',e_2'',e_3',e_3'',e_5'$ are
terminal edges.
\end{enumerate}

By Lemma~\ref{GammaMType52}(b),
either  $w(\Gamma\cap {\rm Int}D_2)=0$ or
 $w(\Gamma\cap {\rm Int}D_3)=0$.
If $w(\Gamma\cap {\rm Int}D_2)=0$,
then by (2) and (4),
we have $e_2'\ni w_4$ and $e_2''\ni w_5$.
Thus $e_2'\cup e_4$ bounds a lens in $D_2$.
This contradicts Lemma~\ref{NoLens}.
If $w(\Gamma\cap {\rm Int}D_3)=0$,
then by (3) and (4),
we have $e_3'\ni w_4$ and $e_3''\cap e_5'\ni w_1$.
Thus $e_5'\cup e_2$ bounds a lens in $D_3$.
This contradicts Lemma~\ref{NoLens}.
Hence Case (i) does not occur.

{\bf Case (ii).}
By (1),
the two terminal edges $e$ and $e'$ are oriented inward 
at $w_1,w_2$,
respectively.
This contradicts Lemma~\ref{Special4Or5angledDiskTwoFeeler}.
Thus Case (ii) does not occur.

{\bf Case (iii).}
We use the notations as shown in Fig.~\ref{Fig27}(b),
where
\begin{enumerate}
\item[(5)] $e_3',e_3'',e_4',e_5''$ are internal edges 
(possibly terminal edges) of label $m+1$ 
oriented outward at $w_3,w_3,w_4,w_5$ in $D_2$,
respectively, and
\end{enumerate}
none of the four edges are middle at $w_3,w_3,w_4,w_5$,
respectively.
Thus by Assumption~\ref{AssumeTerminal},
none of the four edges are
terminal edges.
Hence we have $w(\Gamma\cap{\rm Int}D_2)\ge1$
by IO-Calculation with respect to $\Gamma_{m+1}$ in $D_2$.
Thus by Lemma~\ref{GammaMType52},
the disk $D_2$ contains an oval of label $m$ as shown in 
Fig.~\ref{Fig11},
where $k=m$ and $\varepsilon=+1$.

Let $D$ be the 2-angled disk of $\Gamma_m$ in $D_2$
without feelers.
By IO-Calculation with respect to $\Gamma_{m+1}$ in $Cl(D_2-D)$,
we have $w(\Gamma\cap({\rm Int}D_2-D))\ge1$.
Thus $w(\Gamma\cap{\rm Int}D_2)\ge3$.
This contradicts Lemma~\ref{GammaMType52}(d).
Hence Case (iii) does not occur.

Therefore all of the three cases do not occur.
Hence $\Gamma_m$
does not contain the graph as shown in Fig.~\ref{Fig12}(b).

Similarly we can show that
 $\Gamma_{m+1}$
does not contain the graph as shown in Fig.~\ref{Fig12}(b).
\end{Proof}

\begin{figure}
\centerline{\includegraphics{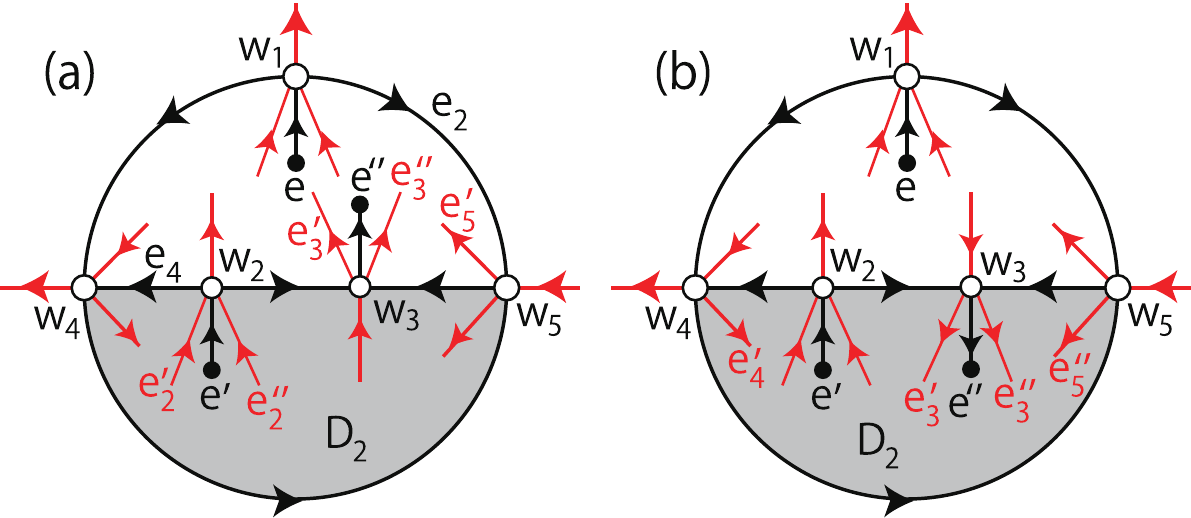}}
\caption{\LABEL{Fig27} 
The gray regions are the disk $D_2$.}
\end{figure}


\section{Case of the graph as shown in Fig.~\ref{Fig12}(c)}

\LABEL{s:TypeC}

In this section,
we shall show that
neither $\Gamma_{m}$ nor $\Gamma_{m+1}$ contains
the graph as shown in Fig.~\ref{Fig12}(c)
for any minimal chart $\Gamma$ of type $(m;7)$.

\begin{lemma}
\LABEL{NoTypeC}
Let $\Gamma$ be a minimal chart of type $(m;7)$.
Then neither $\Gamma_m$ nor $\Gamma_{m+1}$
contains the graph as shown in Fig.~\ref{Fig12}$($c$)$.
\end{lemma}

\begin{Proof}
Suppose that $\Gamma_m$
contains the graph as shown in Fig.~\ref{Fig12}(c),
say $G$.
By Lemma~\ref{OriGammaM5}(b),
if necessary we
move the point at infinity $\infty$, then
the graph $G$ is one of the RO-family of the graph
as shown in Fig.~\ref{Fig13}(c).
Without loss of generality
we can assume that
the graph $G$ is the graph
as shown in Fig.~\ref{Fig13}(c).

We use the notations as shown in Fig.~\ref{Fig14}(b),
where $w_3,w_4,w_5$ are BW-vertices.
Now the graph $G$ divides $S^2$ into three disks.
Let $D_1,D_2,D_3$ be the three disks
which are 4-angled disks
with $\partial D_1\ni w_3,w_4$,
$\partial D_2\ni w_3,w_5$,
$\partial D_3\ni w_4,w_5$.
Let $e,e',e''$ be the terminal edges of label $m$
at $w_3,w_4,w_5$, respectively.
Then $e',e''$ are oriented outward at $w_4,w_5$,
respectively (see Fig.~\ref{Fig13}(c)).
Since $D_3$ is a special 4-angled disk,
by Lemma~\ref{Special4Or5angledDiskTwoFeeler}
we have $e'\not\subset D_3$ or $e''\not\subset D_3$.
Without loss of generality
we can assume that
\begin{enumerate}
\item[(1)] $e'\not\subset D_3$ (i.e. $e'\subset D_1$).
\end{enumerate}
For the edge $e$,
there are two cases:
(i) $e\subset D_1$ (see Fig.~\ref{Fig28}(a)),
(ii) $e\subset D_2$ (see Fig.~\ref{Fig28}(b) and (c)).

{\bf Case (i).}
We use the notations as shown in Fig.~\ref{Fig28}(a),
where
\begin{enumerate}
\item[(2)] $e_1',e_2',e_3',e_3''$ are internal edges 
(possibly terminal edges) of label $m+1$ 
oriented inward at $w_1,w_2,w_3,w_3$ in $D_1$,
respectively, and
\end{enumerate}
none of the four edges are middle at $w_1,w_2,w_3,w_3$,
respectively.
Thus by Assumption~\ref{AssumeTerminal},
none of the four edges are terminal edges.
Hence we have $w(\Gamma\cap{\rm Int}D_1)\ge1$
by IO-Calculation with respect to $\Gamma_{m+1}$ in $D_1$.
Thus by Lemma~\ref{GammaMType52},
the disk $D_1$ contains an oval of label $m$ as shown in 
Fig.~\ref{Fig11},
where $k=m$ and $\varepsilon=+1$.

Let $D$ be the 2-angeled disk of $\Gamma_m$ in $D_1$
without feelers.
By IO-Calculation with respect to $\Gamma_{m+1}$ in $Cl(D_1-D)$,
we have $w(\Gamma\cap({\rm Int}D_1-D))\ge1$.
Thus $w(\Gamma\cap{\rm Int}D_1)\ge3$.
This contradicts Lemma~\ref{GammaMType52}(d).
Hence Case (i) does not occur.

{\bf Case (ii).}
Now $e''\subset D_3$.
Because if $e''\subset D_2$ (see Fig.~\ref{Fig28}(b)),
then we can show $w(\Gamma\cap{\rm Int}D_2)\ge3$ 
by the similar way of the proof of Case (i).
This contradicts Lemma~\ref{GammaMType52}(d).
Hence we have $e''\subset D_3$ (see Fig.~\ref{Fig28}(c)).

We use the notations as shown in Fig.~\ref{Fig28}(c),
where
\begin{enumerate}
\item[(3)] $e_1'',e_2'',e_3',e_3''$ are internal edges 
(possibly terminal edges) of label $m+1$ 
oriented inward at $w_1,w_2,w_3,w_3$ in $D_2$,
respectively, and
\end{enumerate}
none of the four edges are middle at $w_1,w_2,w_3,w_3$,
respectively.
Thus by Assumption~\ref{AssumeTerminal},
none of the four edges are terminal edges.
Hence we can show that $w(\Gamma\cap{\rm Int}D_2)\ge3$
by the similar way of the proof of Case (i).
This contradicts Lemma~\ref{GammaMType52}(d).
Hence Case (ii) does not occur.

Therefore the two cases do not occur.
Hence $\Gamma_m$
does not contain the graph as shown in Fig.~\ref{Fig12}(c).

Similarly we can show that
 $\Gamma_{m+1}$
does not contain the graph as shown in Fig.~\ref{Fig12}(c).
\end{Proof}

\begin{figure}
\centerline{\includegraphics{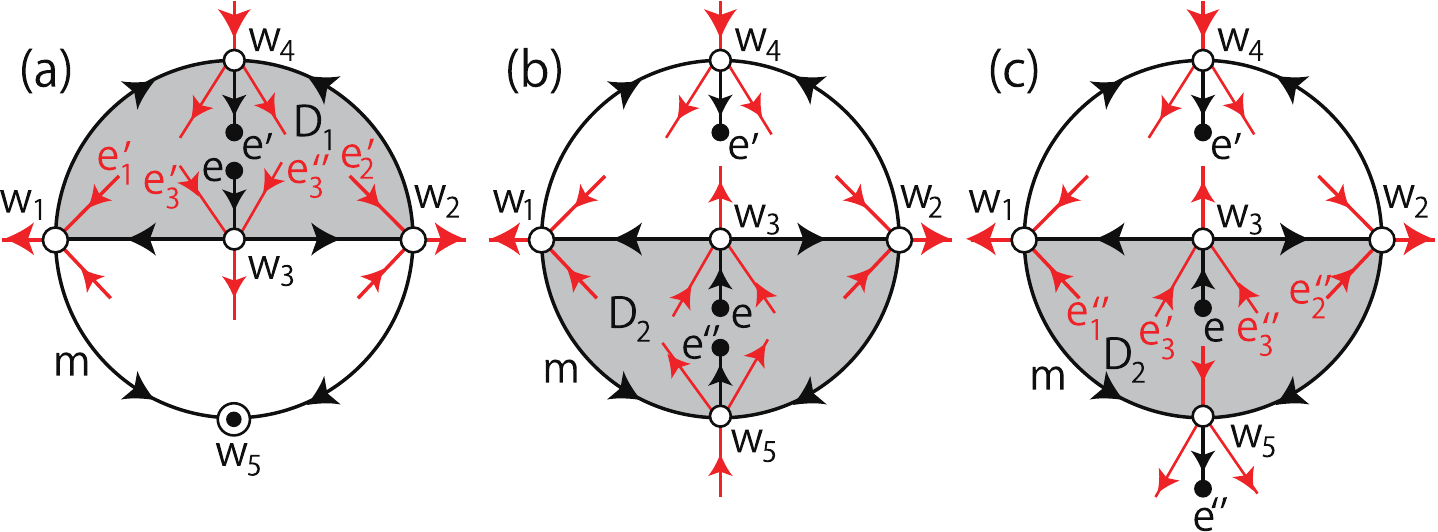}}
\caption{\LABEL{Fig28} 
(a) The gray region is the disk $D_1$.
(b),(c) The gray regions are the disk $D_2$.}
\end{figure}

\begin{proposition}
\LABEL{PropTypeG}
If there exists a minimal chart $\Gamma$ of type $(m;7)$,
then both of $\Gamma_m$ and $\Gamma_{m+1}$
contain  graphs as shown in Fig.~\ref{Fig12}$($g$)$.
\end{proposition}

\begin{Proof}
By Lemma~\ref{GammaMType52},
both of $\Gamma_m$ and $\Gamma_{m+1}$ are of type $(5,2)$.
By Lemma~\ref{LemmaNoLoop} and Lemma~\ref{LemmaWithTerminal}, 
each of $\Gamma_m$ and $\Gamma_{m+1}$ contains
one of nine graphs as shown in Fig.~\ref{Fig12}.
By the eight lemmata
(Lemma~\ref{NoTypeA},
Lemma~\ref{NoTypeD},
Lemma~\ref{NoTypeF},
Lemma~\ref{NoTypeI},
Lemma~\ref{NoTypeH},
Lemma~\ref{NoTypeE},
Lemma~\ref{NoTypeB},
Lemma~\ref{NoTypeC}),
both of $\Gamma_m$ and $\Gamma_{m+1}$
contain graphs as shown in Fig.~\ref{Fig12}$($g$)$.
\end{Proof}


\section{Proof of Main Theorem}
\LABEL{s:Main}

In this section we shall prove Main Theorem(Theorem~\ref{MainTheorem}).

\begin{lemma}
\LABEL{LemmaStep2TypeG}
Let $\Gamma$ be a minimal chart of type $(m;7)$.
Then $\Gamma$ contains one of the RO-family of 
the pseudo chart as
shown in Fig.~\ref{Fig29}$($a$)$.
Moreover the edge $e_3$ is middle at $w_2$,
but not middle at $w_3$.
\end{lemma}

\begin{Proof}
By Proposition~\ref{PropTypeG},
the graph $\Gamma_m$ contains the graph as shown in 
Fig.~\ref{Fig12}(g), say $G$.
By Lemma~\ref{OriGammaM5}(c),
the graph $G$ is one of the RO-family of the graph
as shown in Fig.~\ref{Fig13}(f).
Without loss of generality
we can assume that
the graph $G$ is the graph
as shown in Fig.~\ref{Fig13}(f).

Let $D_1,D_2$ be the 2-angled disk of $\Gamma_m$
and the 3-angled disk of $\Gamma_m$
with $\partial D_1\subset G$,
$\partial D_2\subset G$ and
$D_1\cap D_2=\emptyset$.
Since $D_1$ is a special 2-angled disk with at most one feeler,
the disk $D_1$ has no feelers
by Lemma~\ref{No2angledDiskOneFeeler}.
Since $D_2$ is a special 3-angled disk,
the disk $D_2$ has no feelers
by Lemma~\ref{Special3angledDiskOneOrTwoFeelers}.
Hence $\Gamma$ contains the pseudo chart
as shown in Fig.~\ref{Fig29}(b).

We use the notations as shown in 
Fig.~\ref{Fig29}(b),
where
\begin{enumerate}
\item[(1)] $e_2',e_2'',e_3',e_4',e_4''$ are internal edges 
(possibly terminal edges) of label $m+1$ 
oriented outward at $w_2,w_2,w_3,w_4,w_4$,
respectively.
\end{enumerate}
Moreover, 
none of the five edges are middle at $w_2,w_2,w_3,w_4,w_4$,
respectively.
Thus by Assumption~\ref{AssumeTerminal},
\begin{enumerate}
\item[(2)] none of the five edges $e_2',e_2'',e_3',e_4',e_4''$ are terminal edges.
\end{enumerate}

{\bf Claim.}
$w(\Gamma\cap(S^2-(D_1\cup D_2)))\ge1$.

{\it Proof of Claim.}
Suppose $w(\Gamma\cap(S^2-(D_1\cup D_2)))=0$.
By (2), the edge $e_3'$ is not a terminal edge.
Thus for the edge $e_3'$,
there are four cases:
(i) $e_3'=e_1'$,
(ii) $e_3'=e_1''$,
(iii) $e_3'$ is a loop,
(iv) $e_3'\ni w_5$.

{\bf Case (i).}
The annulus $Cl(S^2-(D_1\cup D_2))$ is 
divided by $e_3\cup e_3'$
into two regions.
Let $E$ be one of the two regions with $E\supset e_2'$.
Then we have $w(\Gamma\cap{\rm Int}E)\ge1$
by IO-Calculation with respect to $\Gamma_{m+1}$
in $E$.
Hence
$$0=w(\Gamma\cap(S^2-(D_1\cup D_2)))\ge w(\Gamma\cap{\rm Int}E)\ge1.$$
This is a contradiction.
Hence Case (i) does not occur.

{\bf Case (ii).}
The condition $e_3'=e_1''$ implies 
$e_1'=e_2'$.
Thus $e_1'\cup e_1$ bounds a lens.
This contradicts Lemma~\ref{NoLens}.
Hence Case (ii) does not occur.

{\bf Case (iii).}
By Lemma~\ref{LemmaNoLoop},
Case (iii) does not occur.

{\bf Case (iv).}
The annulus $Cl(S^2-(D_1\cup D_2))$ is
divided by $e_3'$
into two regions.
Let $E$ be one of the two regions with $E\supset e_4'$.
Then we have $w(\Gamma\cap{\rm Int}E)\ge1$
by IO-Calculation with respect to $\Gamma_{m+1}$
in $E$.
Thus we have the same contradiction
by the similar way of Case (i).
Hence Case (iv) does not occur.

Therefore all the four cases do not occur.
Hence Claim holds. \hfill {$\square$}\vspace{1.5em}

By Claim and Lemma~\ref{GammaMType52}(c),
the region $S^2-(D_1\cup D_2)$ contains an oval of
label $m$ as shown in Fig.~\ref{Fig11},
where $k=m$ and $\varepsilon=+1$.
Moreover by Lemma~\ref{GammaMType52},
we have 
\begin{enumerate}
\item[(3)] $w(\Gamma\cap{\rm Int}D_1)=0$ and $w(\Gamma\cap{\rm Int}D_2)=0$.
\end{enumerate}
Thus by Lemma~\ref{Theorem2AngledDisk},
a regular neighborhood of $D_1$ contains
the pseudo chart as shown in Fig.~\ref{Fig03}(b).
Furthermore,
 by Lemma~\ref{Theorem3AngledDisk},
a regular neighborhood of $D_2$ contains one of 
the RO-family of 
the pseudo chart as shown in Fig.~\ref{Fig09}(b).
Therefore $\Gamma$ contains the pseudo chart as
shown in Fig.~\ref{Fig29}$($a$)$.
\end{Proof}

\begin{figure}
\centerline{\includegraphics{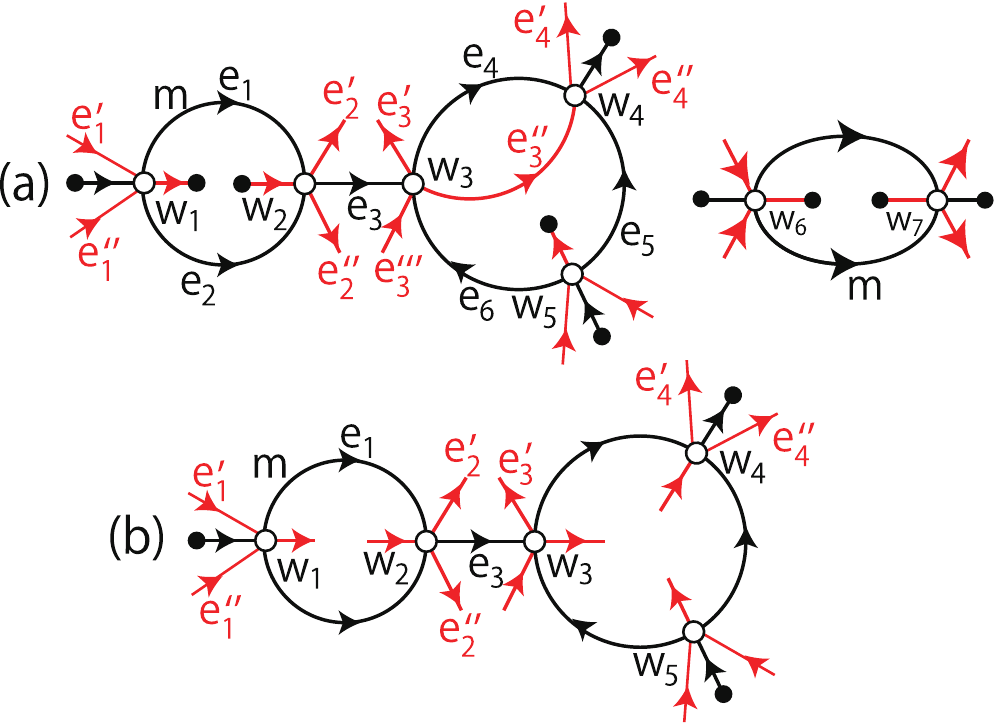}}
\caption{\LABEL{Fig29} 
The graphs as shown in Fig.~\ref{Fig12}(g).}
\end{figure}

{\it Proof of Main Theorem$($Theorem~\ref{MainTheorem}$)$.}
Suppose that there exists a minimal chart $\Gamma$
of type $(m;7)$.
Then by Lemma~\ref{LemmaStep2TypeG},
the chart $\Gamma$ contains one of the RO-family of the
pseudo chart as shown in Fig.~\ref{Fig29}(a),
where 
\begin{enumerate}
\item[(1)] $e_3''$ is an internal edge of label $m+1$
middle at $w_4$ with $\partial e_3''=\{w_3,w_4\}$.
\end{enumerate}

By Proposition~\ref{PropTypeG},
the graph $\Gamma_{m+1}$ contains the graph as shown in
Fig.~\ref{Fig12}(g).
Moreover by Lemma~\ref{GammaMType52},
the graph $\Gamma_{m+1}$ contains an oval of label $m+1$
as shown in Fig.~\ref{Fig11},
where $k=m+1$ and $\varepsilon=-1$.
By the similar way of the proof of Lemma~\ref{LemmaStep2TypeG},
the chart $\Gamma$ contains the graph as shown in 
Fig.~\ref{Fig30}(a), and
\begin{enumerate}
\item[(2)] $e$ is an internal edge of label $m+1$
middle at $v_2$ with $\partial e=\{v_2,v_3\}$,
\item[(3)]
there exists an internal edge $e'$ of label $m$
with $\partial e'=\{v_3,v_4\}$ or $\partial e'=\{v_3,v_5\}$,
\item[(4)] there exists a disk $D$ such that 
$\partial D$ consists of the edge $e'$ and 
an internal edge of label $m+1$ and
${\rm Int}D$  does not intersect any edge at $v_3,v_4$ or $v_5$.
\end{enumerate}

{\bf Claim 1.}
$e_3''=e$,  $w_3=v_3$ and $w_4=v_2$.

{\it Proof of Claim $1$.}
First we show $e_3''=e$.
Since $\Gamma_{m+1}$ contains exactly five terminal edges
(see Fig.~\ref{Fig30}(a)) and
since the five white vertices 
$w_1,w_2,w_5,w_6,w_7$ are contained in terminal edges 
of label $m+1$ (see Fig.~\ref{Fig29}(a)),
 the two white vertices $w_3,w_4$ are not contained in any terminal edge 
of label $m+1$. 
Thus $\{w_3,w_4\}=\{v_2,v_3\}$ and
$e_3''=e$.

Since the edges $e_3''$ and $e$ are middle at $w_4,v_2$,
respectively,
by (1) and (2),
we have $w_4=v_2$.
This implies $w_3=v_3$.
Hence Claim~1 holds.
{\hfill {$\square$}}\vspace{1.5em}

Let $C,C'$ be the simple closed curves in $\Gamma_{m+1}$
with $v_1,v_2\in C$ and $v_3,v_4,v_5\in C'$
(see Fig.~\ref{Fig30}(a)).

{\bf Claim 2.} 
$C=e_4'\cup e_4''\ni w_4$, $e_3'''=e_2''$ and
$C'=e_3'\cup e_3'''\cup e_2'\ni w_2,w_3$.

{\it Proof of Claim~$2$.}
By Claim~1,
we have $w_4=v_2\in C$ and $w_3=v_3\in C'$.
Hence $w_4\in C$ and $w_3\in C'$.
Thus $C=e_4'\cup e_4''$ and $C'\supset e_3'\cup e_3'''$.

Next we shall show $w_2\in C'$.
Since $\partial e'=\{v_3,v_4\}$ or 
$\partial e'=\{v_3,v_5\}$ by (3), 
we have $\partial e'\subset C'$.

Since $e_3''=e$ by Claim 1, we have $e_3=e'$
(see Fig.~\ref{Fig29}(a) and
Fig.~\ref{Fig30}(a)).
Thus $\partial e'=\partial e_3=\{w_2,w_3\}$.
Hence the condition $\partial e'\subset C'$ implies $w_2\in C'$.

Thus $C'\supset e_2'\cup e_2''$.
Hence $e_3'''=e_2'$ or $e_3'''=e_2''$.
Thus by (4), we have $e_3'''=e_2''$ and 
 $C'=e_3'\cup e_3'''\cup e_2'$.
Therefore Claim~2 holds.
{\hfill {$\square$}}\vspace{1.5em}

Now the simple closed curve $C$ contains the two white vertices
$v_1$ and $v_2=w_4$.
For the white vertex $v_1$,
there are three cases:
(i) $v_1=w_1$,
(ii) $v_1=w_5$,
(iii) $v_1=w_6$.

{\bf Case (i).}
Let $D_1,D_2$ be the 2-angled disk of $\Gamma_m$
and the 3-angled disk of $\Gamma_m$
with $w_1,w_2\in \partial D_1$,
$w_3,w_4,w_5\in \partial D_2$ and
$D_1\cap D_2=\emptyset$.
Then the annulus $Cl(S^2-(D_1\cup D_2))$
is divided by $e_4'\cup e_3$
into two disks.
Let $E$ be one of the two disks
with $E\supset e_2'$.
Then $w(\Gamma\cap {\rm Int}E)\ge1$
by IO-Calculation with respect to $\Gamma_{m+1}$
in $E$.
Hence $E$ contains the oval of label $m$.

Let $D'$ be the 2-angled disk of $\Gamma_m$ in $E$.
Then we can show $w(\Gamma\cap(S^2-(D_1\cup D_2)))\ge3$
by using IO-Calculation with respect to $\Gamma_{m+1}$
in ${\rm Int}E-D'$.
This contradicts Lemma~\ref{GammaMType52}(d).
Hence Case (i) does not occur.

{\bf Case (ii).}
The edges $e_4',e_4''$ contain the white vertex $w_5$.
Thus $e_4''\cup e_5$ bounds a lens.
This contradicts Lemma~\ref{NoLens}.
Hence Case (ii) does not occur.

{\bf Case (iii).}
The edges $e_4',e_4''$ contain the white vertex $w_6$.
By Claim~2, we have $e_3'''=e_2''$.
Hence the two edges $e_2',e_3'$ contain  the white vertex $w_1$
(see Fig.~\ref{Fig30}(b)).
Thus $e_2'\cup e_1$ bounds a lens.
This contradicts Lemma~\ref{NoLens}.
Hence Case (iii) does not occur.

Therefore all the three cases do not occur.
Hence there does not exist a minimal chart of type $(m;7)$.
Thus we complete the proof of Main Theorem(Theorem~\ref{MainTheorem}).
{\hfill {$\square$}}

\begin{figure}
\centerline{\includegraphics{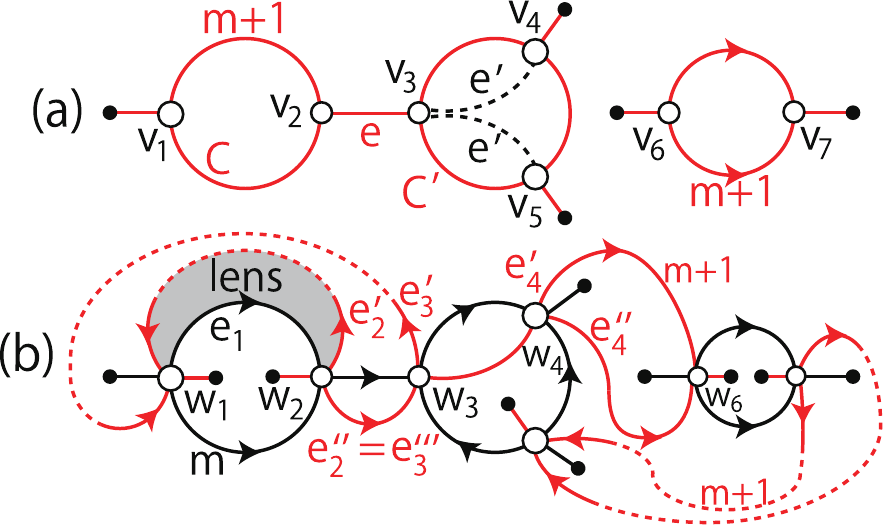}}
\caption{\LABEL{Fig30} 
(a) The graph as shown in Fig.~\ref{Fig12}(g) and an oval.
(b) The gray region is a lens.}
\end{figure}




\vspace{5mm}

\begin{minipage}{65mm}
{Teruo NAGASE
\\
{\small Tokai University \\
4-1-1 Kitakaname, Hiratuka \\
Kanagawa, 259-1292 Japan\\
\\
nagase@keyaki.cc.u-tokai.ac.jp
}}
\end{minipage}
\begin{minipage}{65mm}
{Akiko SHIMA 
\\
{\small Department of Mathematics, 
\\
Tokai University
\\
4-1-1 Kitakaname, Hiratuka \\
Kanagawa, 259-1292 Japan\\
shima@keyaki.cc.u-tokai.ac.jp
}}
\end{minipage}


\vspace{0.7cm}

{\bf List of terminologies}\vspace{2mm}\\
{\small $
\begin{array}{ll||}
\text{$k$-angled disk} & p5 \\
\text{BW-vertex} & p6 \\
\text{C-move equivalent} & p3 \\
\text{chart} & p2 \\
\text{complexity $(w(\Gamma),-f(\Gamma))$} & p3 \\
\text{feeler} & p5 \\
\text{free edge} & p3 \\
\text{hoop} & p4 \\
\text{internal edge} & p6 \\
\text{inward} & p3 \\
\text{inward arc} & p18 \\
\text{IO-Calculation} & p18 \\
\text{keeping $X$ fixed} & p9 \\
\text{lens} & p21 \\
\text{local complexity $\ell c(D;\Gamma)$} & p9 \\
\text{locally minimal} & p9 \\
\text{loop} & p10 \\
\end{array}
~~
\begin{array}{ll}
\text{middle arc} & p3 \\
\text{middle at $v$} & p3 \\
\text{minimal chart} & p3 \\
\text{outward} & p3 \\
\text{outward arc} & p18 \\
\text{oval} & p6 \\
\text{point at infinity $\infty$} & p4 \\
\text{pseudo chart} & p5 \\
\text{ring} & p4 \\
\text{RO-family} & p9 \\
\text{simple hoop} & p4 \\
\text{skew $\theta$-curve} & p6 \\
\text{special $k$-angled disk} & p9 \\
\text{terminal edge} & p4 \\
\text{type $(m;n_1,n_2,\cdots,n_k)$ for a chart} & p1 \\
\text{type $(n_1,n_2,\cdots,n_k)$ for $\Gamma_m$} & p10 \\
\text{$\theta$-curve} & p6 \\
\end{array}
$}

\vspace{0.5cm}

{\bf List of notations}\vspace{2mm}\\
{\small $
\begin{array}{ll}
\text{$\Gamma_m$} & p1 \\
\text{$w(\Gamma)$} & p3 \\
\text{$f(\Gamma)$} & p3 \\
\text{${\rm Int}X$} & p5 \\
\text{$\partial X$} & p5 \\
\text{$Cl(X)$} & p5 \\
\text{$\partial \alpha$} & p5 \\
\text{$w(X)$} & p5 \\
\text{$c(X)$} & p9 \\
\text{$\ell c(D;\Gamma)$} & p9\\
\end{array}
$
}

\end{document}